%% file: cmpDimMZV_Ar_v2.tex
\documentclass[11pt]{article}
	\newcommand{\mainTitle}{Computations about formal multiple zeta spaces defined by binary extended double shuffle relations}	
	\newcommand{\authorName}{Tomoya Machide}	
	\newcommand{\organizationNameFst}{National Institute of Informatics}
	\newcommand{\placeAddressFst}{2-1-2 Hitotsubashi, Chiyoda-ku, Tokyo 101-8430, Japan}
	\newcommand{\emailAddressFst}{machide@nii.ac.jp}

	\newcommand{\MSCname}{11M32 (Primary); 15A03,68W30 (Secondary)} 
	\newcommand{\keyWord}{multiple zeta value, graded vector space, dimension calculation} 

	\newcommand{\siteMyCode}{https://github.com/machide-tomoyan/BMZS-calculator}	

\usepackage{graphicx}
\usepackage{geometry}               
\usepackage{ascmac}
\usepackage{amsmath}
\usepackage{amssymb}
\usepackage{amsthm}
\usepackage[bookmarkstype=toc,colorlinks=true,urlcolor=black,linkcolor=black,citecolor=black,linktocpage=true,bookmarks=false]{hyperref}
\usepackage{arydshln}
\usepackage{multicol}
\usepackage{multirow}

	\DeclareMathOperator*{\OPlus}{\bigoplus}
	\DeclareMathOperator*{\LAnd}{\land}	\DeclareMathOperator*{\LAND}{\bigwedge}	
	\DeclareMathOperator*{\LOr}{\lor}		\DeclareMathOperator*{\LOR}{\bigvee}		
	\DeclareMathOperator*{\dmoPR}{\Pr}

	\newcommand{\RCTbTmp}{{\ooalign{\hfill${\ooalign{\hfill${\ooalign{\hfill$\!\mid$\hfill\crcr\hfill$\!\mid$\hfill}}$\hfill\crcr\hfill$\!\!\mid$\hfill}}$\hfill\crcr\hfill\hspace{-5pt}$\mid$\hfill}}}
	\newcommand{\RCT}[1][n]{\ifx n#1\!{\ooalign{\hfill$\lceil$\hfill\crcr\hfill$\>\>\rfloor$\hfill}}\else\rotatebox[origin=c]{90}{{\ooalign{\hfill\rotatebox[origin=c]{180}{$\neg$}\hfill\crcr\hfill$\neg$\hfill}}}\fi}	
	\newcommand{\bRCT}[1][n]{\ifx n#1{\ooalign{\hfill$\RCTbTmp$\hfill\crcr\hfill$\hspace{1.5pt}\RCTbTmp$\hfill}}
						       \else\rotatebox[origin=c]{90}{{\ooalign{\hfill\rotatebox[origin=c]{180}{$\neg$}\hfill\crcr\hfill$\neg$\hfill\crcr\hfill{\small\!$\bullet$}\hfill\crcr\hfill{\small\,$\bullet$}\hfill}}}\fi}	
	\newcommand{\nbk}[3]{#1#3#2}		
	\newcommand{\bgbk}[3]{\bigl{#1}#3\bigr{#2}}	
	\newcommand{\Bgbk}[3]{\Bigl{#1}#3\Bigr{#2}}			
	\newcommand{\bggbk}[3]{\biggl{#1}#3\biggr{#2}}			
	\newcommand{\Bggbk}[3]{\Biggl{#1}#3\Biggr{#2}}
	\newcommand{\autobk}[3]{\left#1#3\right#2}
	\newcommand{\nbkD}[5]{#1#2#5#3#4}		
	\newcommand{\bgbkD}[5]{\bigl{#1}\bigl{#2}#5\bigr{#3}\bigr{#4}}	
	\newcommand{\BgbkD}[5]{\Bigl{#1}\Bigl{#2}#5\Bigr{#3}\Bigr{#4}}	
	\newcommand{\bggbkD}[5]{\biggl{#1}\biggl{#2}#5\biggr{#3}\biggr{#4}}	
	\newcommand{\BggbkD}[5]{\Biggl{#1}\Biggl{#2}#5\Biggr{#3}\Biggr{#4}}	
	\newcommand{\autobkD}[5]{\left#1\left#2#5\right#3\right#4}	
	\newcommand{\mcbk}[4][?]{\ifx n#1\nbk{#2}{#3}{#4}\else\ifx b#1\bgbk{#2}{#3}{#4}\else\ifx B#1\Bgbk{#2}{#3}{#4}\else\ifx g#1\bggbk{#2}{#3}{#4}\else\ifx G#1\Bggbk{#2}{#3}{#4}\else\ifx a#1\autobk{#2}{#3}{#4}\else\ifx !#1{#4}\else#4\fi\fi\fi\fi\fi\fi\fi}
	\newcommand{\mcbkD}[4][?]{\ifx n#1\nbkD{#2}{#2}{#3}{#3}{#4}\else\ifx b#1\bgbkD{#2}{#2}{#3}{#3}{#4}\else\ifx B#1\BgbkD{#2}{#2}{#3}{#3}{#4}\else\ifx g#1\bggbkD{#2}{#2}{#3}{#3}{#4}\else\ifx G#1\BggbkD{#2}{#2}{#3}{#3}{#4}\else\ifx a#1\autobkD{#2}{#2}{#3}{#3}{#4}\else\ifx !#1{#4}\else#4\fi\fi\fi\fi\fi\fi\fi}
	\newcommand{\nsgsb}[1]{#1}		
	\newcommand{\bgsgsb}[1]{\big{#1}}	
	\newcommand{\Bgsgsb}[1]{\Big{#1}}			
	\newcommand{\bggsgsb}[1]{\bigg{#1}}			
	\newcommand{\Bggsgsb}[1]{\Bigg{#1}}
	\newcommand{\mcsgsb}[2][?]{\ifx n#1\nsgsb{#2}\else\ifx b#1\bgsgsb{#2}\else\ifx B#1\Bgsgsb{#2}\else\ifx g#1\bggsgsb{#2}\else\ifx G#1\Bggsgsb{#2}\else#2\fi\fi\fi\fi\fi}
	\newcommand{\myEqSpace}{\,} 	\newlength{\myEqSpaceLen} 	
	\newcommand{\eLt}[1]{\widehat{#1}}
	\newcommand{\rLt}[1]{\overline{#1}}
	\newcommand{\bLt}[1]{\underline{#1}}

	\newcommand{\bkR}[2][n]{\mcbk[#1]{(}{)}{#2}}						
	\newcommand{\bkS}[2][n]{\mcbk[#1]{[}{]}{#2}}						
	\newcommand{\bkB}[2][n]{\mcbk[#1]{\{}{\}}{#2}}						
	\newcommand{\bkA}[2][n]{\mcbk[#1]{\langle}{\rangle}{#2}}				
	\newcommand{\bkF}[2][n]{\mcbk[#1]{\lfloor}{\rfloor}{#2}}				
	\newcommand{\bkAll}[4][n]{\mcbk[#1]{#2}{#3}{#4}}
		
	\newcommand{\nFc}[3][n]{#2\bkR[#1]{#3}}					
				
	\newcommand{\idFc}[4][n]{\id{#2}{#3}\bkR[#1]{#4}}			
	\newcommand{\pwFc}[4][n]{\pw{#2}{#3}\bkR[#1]{#4}}			
	\newcommand{\ipFc}[5][n]{\ip{#2}{#3}{#4}\bkR[#1]{#5}}		
		\newcommand{\Fc}{\nFc}		
			
	\newcommand{\del}{\delta}
	
	\newcommand{\gam}{\gamma} 
	\newcommand{\Gam}{\Gamma} 
	
	\newcommand{\ep}{\varepsilon}

	\newcommand{\modd}[1][\ ]{\mathrm{mod}#1}
	
	\newcommand{\bfAlp}[1][s]{\ifx s#1{\boldsymbol\alpha}\else{\boldsymbol??none??}\fi}		
	\newcommand{\bfBeta}[1][s]{\ifx s#1{\boldsymbol\beta}\else{\boldsymbol??none??}\fi}		
	\newcommand{\bfDelta}[1][s]{\ifx s#1{\boldsymbol\delta}\else{\boldsymbol??none??}\fi}		
	\newcommand{\bfGam}[1][s]{\ifx s#1{\boldsymbol\gam}\else{\boldsymbol\Gam}\fi}			
	\newcommand{\bfA}[1][s]{\ifx s#1{\bf a}\else{\bf A}\fi}
	\newcommand{\bfB}[1][s]{\ifx s#1{\bf b}\else{\bf B}\fi}
	\newcommand{\bfC}[1][s]{\ifx s#1{\bf c}\else{\bf C}\fi}
	\newcommand{\bfD}[1][s]{\ifx s#1{\bf d}\else{\bf D}\fi}
	\newcommand{\bfE}[1][s]{\ifx s#1{\bf e}\else{\bf E}\fi}
	\newcommand{\bfH}[1][s]{\ifx s#1{\bf h}\else{\bf H}\fi}
	\newcommand{\bfI}[1][s]{\ifx s#1{\bf i}\else{\bf I}\fi}
	\newcommand{\bfJ}[1][s]{\ifx s#1{\bf j}\else{\bf J}\fi}
	\newcommand{\bfK}[1][s]{\ifx s#1{\bf k}\else{\bf K}\fi}
	\newcommand{\bfL}[1][s]{\ifx s#1{\bf l}\else{\bf L}\fi}
	\newcommand{\bfM}[1][s]{\ifx s#1{\bf m}\else{\bf M}\fi}
	\newcommand{\bfN}[1][s]{\ifx s#1{\bf n}\else{\bf N}\fi}
	\newcommand{\bfP}[1][s]{\ifx s#1{\bf p}\else{\bf P}\fi}
	\newcommand{\bfQ}[1][s]{\ifx s#1{\bf q}\else{\bf Q}\fi}
	\newcommand{\bfR}[1][s]{\ifx s#1{\bf r}\else{\bf R}\fi}
	\newcommand{\bfS}[1][s]{\ifx s#1{\bf s}\else{\bf S}\fi}
	\newcommand{\bfT}[1][s]{\ifx s#1{\bf t}\else{\bf T}\fi}
	\newcommand{\bfU}[1][s]{\ifx s#1{\bf u}\else{\bf U}\fi}
	\newcommand{\bfV}[1][s]{\ifx s#1{\bf v}\else{\bf V}\fi}
	\newcommand{\bfW}[1][s]{\ifx s#1{\bf w}\else{\bf W}\fi}
	\newcommand{\bfX}[1][s]{\ifx s#1{\bf x}\else{\bf X}\fi}		
	\newcommand{\bfY}[1][s]{\ifx s#1{\bf y}\else{\bf Y}\fi}		
	\newcommand{\bfZ}[1][s]{\ifx s#1{\bf z}\else{\bf Z}\fi}		
	\newcommand{\bfZero}[1][?]{{\bf 0}}
	\newcommand{\bfEp}[1][s]{\ifx s#1{\boldsymbol\ep}\else{\boldsymbol{\mathcal{E}}}\fi}
	\newcommand{\mo}{(-1)}

	\newcommand{\larr}[1][n]{\ifx#1d\dashleftarrow\else\ifx#1c\curvearrowleft\else\leftarrow\fi\fi}
	\newcommand{\rarr}[1][n]{\ifx#1d\dashrightarrow\else\ifx#1c\curvearrowright\else\rightarrow\fi\fi}
	\newcommand{\mVert}[1][n]{{\,\mcsgsb[#1]{\vert}\,}}		

	\newcommand{\SetO}[2][n]{\bkB[#1]{#2}}
	\newcommand{\SetT}[3][n]{\bkB[#1]{#2\mVert[#1]#3}}
		
	\newcommand{\SpO}[2][n]{\bkA[#1]{#2}}
	\newcommand{\SpT}[3][n]{\bkA[#1]{#2\mVert#3}}
		
	\newcommand{\setN}{\mathbb{N}}
	\newcommand{\setZ}{\mathbb{Z}} 
	\newcommand{\setQ}{\mathbb{Q}}	
	\newcommand{\setR}{\mathbb{R}}	
	
	\newcommand{\setF}[1][?]{\ifx ?#1\mathbb{F}\else\mathbb{F}_{#1}\fi}
	
	\newcommand{\matI}[1][?]{\ifx #1?I\else I_{#1}\fi}	
	
	\newcommand{\gpKleinF}[1][?]{V}
		
	\newcommand{\gpu}[1][?]{\ifx?#1e\else e_{#1}\fi}		

	\newcommand{\vPack}[1][10]{\vspace{-#1pt}}
	\newcommand{\vWiden}[1][10]{\vspace{#1pt}}
	\newcommand{\lnA}[1][]{&  &}
	\newcommand{\lnP}[2][1]{\ifx1#1\myEqSpace#2\myEqSpace\else\myEqSpace\myEqSpace#2\myEqSpace\myEqSpace\fi}
	\newcommand{\lnAP}[2][]{& #2 &}

		\newcommand{\slnAH}[1][?]{\\}
		
	\newcommand{\refEq}[1]{(\ref{#1})}	
		
	\newcommand{\pcstSpForRefThm}{\;}		
	\newcommand{\refHL}[2]{#1\pcstSpForRefThm\ref{#2}}		
	\newcommand{\refHLm}[3][?]{\ifx?#1#2\pcstSpForRefThm#3\else#2#3\fi}
	\newcommand{\refThm}[2][?]{\ifx?#1\refHL{Theorem}{#2}\else\ifx s#1\refHL{Theorems}{#2}\else{[argument error]}\fi\fi}
	\newcommand{\refProp}[2][?]{\ifx?#1\refHL{Proposition}{#2}\else\ifx s#1\refHL{Propositions}{#2}\else{[argument error]}\fi\fi}
	\newcommand{\refLem}[2][?]{\ifx?#1\refHL{Lemma}{#2}\else\ifx s#1\refHL{Lemmas}{#2}\else{[argument error]}\fi\fi}
	\newcommand{\refCor}[2][?]{\ifx?#1\refHL{Corollary}{#2}\else\ifx s#1\refHL{Corollaries}{#2}\else{[argument error]}\fi\fi}
	\newcommand{\refDef}[2][?]{\ifx?#1\refHL{Definition}{#2}\else\ifx s#1\refHL{Definitions}{#2}\else{[argument error]}\fi\fi}
	\newcommand{\refRem}[2][?]{\ifx?#1\refHL{Remark}{#2}\else\ifx s#1\refHL{Remarks}{#2}\else{[argument error]}\fi\fi}
	\newcommand{\refTab}[2][?]{\ifx?#1\refHL{Table}{#2}\else\ifx s#1\refHL{Tables}{#2}\else{[argument error]}\fi\fi}
	\newcommand{\refFig}[2][?]{\ifx?#1\refHL{Figure}{#2}\else\ifx s#1\refHL{Figures}{#2}\else{[argument error]}\fi\fi}
	\newcommand{\refExp}[2][?]{\ifx?#1\refHL{Experiment}{#2}\else\ifx s#1\refHL{Experiments}{#2}\else{[argument error]}\fi\fi}
	\newcommand{\refAlg}[2][?]{\ifx?#1\refHL{Algorithm}{#2}\else\ifx s#1\refHL{Algorithms}{#2}\else{[argument error]}\fi\fi}
	\newcommand{\refPrc}[2][?]{\ifx?#1\refHL{Process}{#2}\else\ifx s#1\refHL{Processes}{#2}\else{[argument error]}\fi\fi}
	\newcommand{\refSec}[2][?]{\ifx?#1\refHL{Section}{#2}\else\ifx s#1\refHL{Sections}{#2}\else{[argument error]}\fi\fi}
	\newcommand{\refApp}[2][?]{\ifx?#1\refHL{Appendix}{#2}\else\ifx s#1\refHL{Appendices}{#2}\else{[argument error]}\fi\fi}
	\newcommand{\refPrb}[2][?]{\ifx?#1\refHL{Problem}{#2}\else\ifx s#1\refHL{Problems}{#2}\else{[argument error]}\fi\fi}
		\newcommand{\refSect}{\refSec}
		\newcommand{\refPrp}{\refProp}
	\newcommand{\refThmA}[2][?]{\ifx?#1\refHLm{Theorem}{#2}\else \refHLm{Theorems}{#2}\fi}
	\newcommand{\refPropA}[2][?]{\ifx?#1 \refHLm[#1]{Proposition}{#2}\else \refHLm{Propositions}{#2}\fi}
	\newcommand{\refLemA}[2][?]{\ifx?#1\refHLm{Lemma}{#2}\else \refHLm{Lemmas}{#2}\fi}
	\newcommand{\refCorA}[2][?]{\ifx?#1\refHLm{Corollary}{#2}\else \refHLm{Corollaries}{#2}\fi}
	\newcommand{\refDefA}[2][?]{\ifx?#1\refHLm{Definition}{#2}\else \refHLm{Definitionss}{#2}\fi}
	\newcommand{\refRemA}[2][?]{\ifx?#1\refHLm{Remark}{#2}\else \refHLm{Remarks}{#2}\fi}
	\newcommand{\refExpA}[2][?]{\ifx?#1\refHLm{Experiment}{#2}\else \refHLm{Experiments}{#2}\fi}
	\newcommand{\refAlgA}[2][?]{\ifx?#1\refHLm{Algorithm}{#2}\else \refHLm{Algorithms}{#2}\fi}
	\newcommand{\refPrcA}[2][?]{\ifx?#1\refHLm{Process}{#2}\else \refHLm{Processes}{#2}\fi}
	\newcommand{\refSecA}[2][?]{\ifx?#1\refHLm{Section}{#2}\else \refHLm{Sections}{#2}\fi}
	\newcommand{\refAppA}[2][?]{\ifx?#1\refHLm{Appendix}{#2}\else \refHLm{Appendices}{#2}\fi}
	\newcommand{\refPrbA}[2][?]{\ifx?#1\refHLm{Problem}{#2}\else \refHLm{Problems}{#2}\fi}

	\newcommand{\frc}[3][?]{\ifx s#1#3/#2\else\ifx b#1(#3)/#2\else\ifx d#1\dfrac{#3}{#2}\else\ifx t#1\tfrac{#3}{#2}\else\frac{#3}{#2}\fi\fi\fi\fi}
	\newcommand{\bnm}[3][?]{\binom{#3}{#2}}

	\newcommand{\pw}[3][?]{\ifx!#3{#2}^{#3}\else#2^{#3}\fi}
	\newcommand{\id}[3][?]{#2_{#3}}
	\newcommand{\ip}[4][?]{{#2}_{#3}^{#4}}
	\newcommand{\pwR}[3][a]{\ifx!#1{\bkR[#1]{#2}}^{#3}\else\bkR[#1]{#2}^{#3}\fi}
	\newcommand{\pwB}[3][a]{\ifx!#1{\bkB[#1]{#2}}^{#3}\else\bkB[#1]{#2}^{#3}\fi}
	\newcommand{\pwS}[3][a]{\ifx!#1{\bkS[#1]{#2}}^{#3}\else\bkS[#1]{#2}^{#3}\fi}

	\newcommand{\nIntO}[2][?]{\ifx l#1\int\limits_{#2}\else\ifx t#1{\textstyle\int\limits_{#2}}\else\int_{#2}\fi\fi}
	\newcommand{\nIntT}[3][?]{\ifx l#1\int\limits_{#2}^{#3}\else\if t#1{\textstyle\int\limits_{#2}^{#3}}\else\int_{#2}^{#3}\fi\fi}	
	\newcommand{\nIntN}[1][?]{\ifx l#1\int\limits\else\ifx t#1{\textstyle\int\limits}\else\int\fi\fi}

	\newcommand{\tpT}[3][a]{ {#2}\atop \bkR[#1]{#3} }

	\newcommand{\nSmO}[2][?]{\ifx l#1\sum\limits_{#2}\else\ifx t#1{\textstyle\sum\limits_{#2}}\else\sum_{#2}\fi\fi}
	\newcommand{\nSmT}[3][?]{\ifx l#1\sum\limits_{#2}^{#3}\else\if t#1{\textstyle\sum\limits_{#2}^{#3}}\else\sum_{#2}^{#3}\fi\fi}	
	\newcommand{\nSmN}[1][?]{\ifx l#1\sum\limits\else\ifx t#1{\textstyle\sum\limits}\else\sum\fi\fi}
	\newcommand{\pSm}[2][?]{\ifx t#1 \sum_{#2}^{\prime} \else \sideset{}{^\prime}\sum_{#2} \fi}
	\newcommand{\pSmT}[3][?]{\ifx t#1 \sum_{#2}^{\prime#3} \else \sideset{}{^\prime}\sum_{#2}^{#3} \fi}	
	\newcommand{\pSmN}[1][?]{\ifx t#1 \sum^{\prime} \else \sideset{}{^\prime}\sum \fi}
	\newcommand{\dSm}[2][?]{\ifx t#1 \sum_{#2}^{\dagger} \else \sideset{}{^\dagger}\sum_{#2} \fi}
	\newcommand{\dSmT}[3][?]{\ifx t#1 \sum_{#2}^{\dagger#3} \else \sideset{}{^\dagger}\sum_{#2}^{#3} \fi}	
	\newcommand{\dSmN}[1][?]{\ifx t#1 \sum^{\dagger} \else \sideset{}{^\dagger}\sum \fi}
	\newcommand{\tpTSm}[3][?]{\nSmO[#1]{\tpT{#2}{#3}}}

		\newcommand{\Sm}{\nSmO}		\newcommand{\SmO}{\nSmO}			\newcommand{\SmT}{\nSmT}			
		\newcommand{\tpSm}{\tpTSm}
		
	\newcommand{\nPd}[2][?]{\ifx l#1 \prod\limits_{#2} \else \prod_{#2} \fi}
	\newcommand{\nPdT}[3][?]{\ifx l#1 \prod\limits_{#2}^{#3} \else \prod_{#2}^{#3} \fi}

	\newcommand{\nOPs}[2][?]{\ifx l#1 \OPlus\limits_{#2} \else \OPlus_{#2} \fi}
	\newcommand{\nOPsT}[3][?]{\ifx l#1 \OPlus\limits_{#2}^{#3} \else \OPlus_{#2}^{#3} \fi}	
	\newcommand{\pOPs}[2][?]{\ifx t#1 \OPlus_{#2}^{\prime} \else \sideset{}{^\prime}\OPlus_{#2} \fi}
	\newcommand{\pOPsT}[3][?]{\ifx t#1 \OPlus_{#2}^{\prime#3} \else \sideset{}{^\prime}\OPlus_{#2}^{#3} \fi}

		\newcommand{\OPs}{\nOPs}

	\newcommand{\nAD}[2][?]{\ifx l#1 \LAND\limits_{#2} \else \LAND_{#2} \fi}
	\newcommand{\nADT}[3][?]{\ifx l#1 \LAND\limits_{#2}^{#3} \else \LAND_{#2}^{#3} \fi}	
	\newcommand{\pAD}[2][?]{\ifx t#1 \LAND_{#2}^{\prime} \else \sideset{}{^\prime}\LAND_{#2} \fi}
	\newcommand{\pADT}[3][?]{\ifx t#1 \LAND_{#2}^{\prime#3} \else \sideset{}{^\prime}\LAND_{#2}^{#3} \fi}

	\newcommand{\nAd}[2][?]{\ifx l#1 \LAnd\limits_{#2} \else \LAnd_{#2} \fi}
	\newcommand{\nAdT}[3][?]{\ifx l#1 \LAnd\limits_{#2}^{#3} \else \LAnd_{#2}^{#3} \fi}	
	\newcommand{\pAd}[2][?]{\ifx t#1 \LAnd_{#2}^{\prime} \else \sideset{}{^\prime}\LAnd_{#2} \fi}
	\newcommand{\pAdT}[3][?]{\ifx t#1 \LAnd_{#2}^{\prime#3} \else \sideset{}{^\prime}\LAnd_{#2}^{#3} \fi}

	\newcommand{\nOR}[2][?]{\ifx l#1 \LOR\limits_{#2} \else \LOR_{#2} \fi}
	\newcommand{\nORT}[3][?]{\ifx l#1 \LOR\limits_{#2}^{#3} \else \LOR_{#2}^{#3} \fi}	
	\newcommand{\pOR}[2][?]{\ifx t#1 \LOR_{#2}^{\prime} \else \sideset{}{^\prime}\LOR_{#2} \fi}
	\newcommand{\pORT}[3][?]{\ifx t#1 \LOR_{#2}^{\prime#3} \else \sideset{}{^\prime}\LOR_{#2}^{#3} \fi}

	\newcommand{\nOr}[2][?]{\ifx l#1 \LOr\limits_{#2} \else \LOr_{#2} \fi}
	\newcommand{\nOrT}[3][?]{\ifx l#1 \LOr\limits_{#2}^{#3} \else \LOr_{#2}^{#3} \fi}	
	\newcommand{\pOr}[2][?]{\ifx t#1 \LOr_{#2}^{\prime} \else \sideset{}{^\prime}\LOr_{#2} \fi}
	\newcommand{\pOrT}[3][?]{\ifx t#1 \LOr_{#2}^{\prime#3} \else \sideset{}{^\prime}\LOr_{#2}^{#3} \fi}

	\newcommand{\nPR}[2][?]{\ifx l#1 \dmoPR\limits_{#2} \else \dmoPR_{#2} \fi}
	\newcommand{\nPRT}[3][?]{\ifx l#1 \dmoPR\limits_{#2}^{#3} \else \dmoPR_{#2}^{#3} \fi}	
	\newcommand{\pPR}[2][?]{\ifx t#1 \dmoPR_{#2}^{\prime} \else \sideset{}{^\prime}\dmoPR_{#2} \fi}
	\newcommand{\pPRT}[3][?]{\ifx t#1 \dmoPR_{#2}^{\prime#3} \else \sideset{}{^\prime}\dmoPR_{#2}^{#3} \fi}

	\newcommand{\nIsO}[2][?]{\ifx l#1 \bigcap\limits_{#2}\else\ifx b#1 \bigcap_{#2}\else{\textstyle\bigcap\limits_{#2}}\fi\fi}
	\newcommand{\nIsT}[3][?]{\ifx l#1 \bigcap\limits_{#2}^{#3}\else\ifx b#1 \bigcap_{#2}^{#3}\else{\textstyle\bigcap\limits_{#2}^{#3}}\fi\fi}	
	\newcommand{\pIs}[2][?]{\ifx t#1 \bigcap_{#2}^{\prime} \else \sideset{}{^\prime}\bigcap_{#2} \fi}
	\newcommand{\pIsT}[3][?]{\ifx t#1 \bigcap_{#2}^{\prime#3} \else \sideset{}{^\prime}\bigcap_{#2}^{#3} \fi}

	\newcommand{\nUnO}[2][l]{\ifx L#1 \bigcup\limits_{#2}\else\ifx b#1 \bigcup_{#2}\else\ifx l#1{\textstyle\bigcup\limits_{#2}}\else{\textstyle\bigcup_{#2}}\fi\fi\fi}
	\newcommand{\nUnT}[3][l]{\ifx L#1 \bigcup\limits_{#2}^{#3}\else\ifx b#1 \bigcup_{#2}^{#3}\else\ifx l#1{\textstyle\bigcup\limits_{#2}^{#3}}\else{\textstyle\bigcup_{#2}^{#3}}\fi\fi\fi}	
	\newcommand{\pUn}[2][l]{\ifx t#1 \bigcup_{#2}^{\prime} \else \sideset{}{^\prime}\bigcup_{#2} \fi}
	\newcommand{\pUnT}[3][l]{\ifx t#1 \bigcup_{#2}^{\prime#3} \else \sideset{}{^\prime}\bigcup_{#2}^{#3} \fi}	
	\newcommand{\dUn}[2][l]{\ifx L#1 \bigsqcup\limits_{#2}\else\ifx b#1 \bigsqcup_{#2}\else\ifx l#1{\textstyle\bigsqcup\limits_{#2}}\else{\textstyle\bigsqcup{#2}}\fi\fi\fi}
	\newcommand{\dUnT}[3][l]{\ifx L#1 \bigsqcup\limits_{#2}^{#3}\else\ifx b#1 \bigsqcup_{#2}^{#3}\else\ifx l#1{\textstyle\bigsqcup\limits_{#2}^{#3}}\else{\textstyle\bigsqcup{#2}^{#3}}\fi\fi\fi}

		\newcommand{\Un}{\nUnO}
		\newcommand{\UnO}{\nUnO}
		\newcommand{\UnT}{\nUnT}

	\newcommand{\nLm}[2][?]{\ifx l#1 \lim\limits_{#2} \else \lim_{#2} \fi}
	\newcommand{\iLm}[2][?]{\ifx l#1 \liminf\limits_{#2} \else \liminf_{#2} \fi}
	\newcommand{\sLm}[2][?]{\ifx l#1 \limsup\limits_{#2} \else \limsup_{#2} \fi}
		
	\newcommand{\nMaxT}[3][a]{ \max\bkB[#1]{#2 \mVert[a] #3} }

		\newcommand{\MaxT}{\nMaxT}

	\newcommand{\nMin}[2][a]{ \min\bkB[#1]{#2} }
	\newcommand{\nMinT}[3][a]{ \min\bkB[#1]{#2 \mVert[a] #3} }

		\newcommand{\Min}{\nMin}		
		\newcommand{\MinT}{\nMinT}

	\newcommand{\dimN}{\dim}
	\newcommand{\dimO}[2][?]{\Fc[#1]{\dimN}{#2}}
	\newcommand{\dimT}[3][?]{\idFc[#1]{\dimN}{#2}{#3}}
	\newcommand{\rnkN}{\mathrm{rank}}
	
	\newcommand{\rnkT}[3][?]{\idFc[#1]{\rnkN}{#2}{\,#3}}

	\newcommand{\glcondEnvLineHead}[1]{ \ifx*#1 \begin{eqnarray*} \else \begin{eqnarray}  \label{#1} \fi }
	\newcommand{\glcondEnvLineTail}[1]{ \ifx*#1 \end{eqnarray*} \else \end{eqnarray} \fi }
	\newcommand{\glcondDis}[1]{\ifx d#1 \displaystyle \fi}
	\newcommand{\glcmdEqShift}{\hspace{-20pt}}
	\newcommand{\glcmdHLineCWiden}{\rule{0cm}{15pt}}	\newcommand{\glcdH}{\glcmdHLineCWiden}
	\newcommand{\lccondPar}[1]{\ifx#1p \\ \fi}

		\newcommand{\envMO}[2][*]{$\ifx d#1 \displaystyle \fi#2$}
		\newcommand{\envMT}[3][*]{$\ifx d#1 \displaystyle \fi#2=#3$}
		\newcommand{\envMTDef}[3][*]{$\ifx d#1 \displaystyle \fi#2:=#3$}
		\newcommand{\envMTPt}[4][*]{$\ifx d#1 \displaystyle \fi#3#2#4$}
			\newcommand{\envM}{\envMT}
			
			\newcommand{\envMPt}{\envMTPt}
		\newcommand{\envMTh}[4][*]{$\ifx d#1 \displaystyle \fi#2=#3=#4$}
		\newcommand{\envMThDef}[4][*]{$\ifx d#1 \displaystyle \fi#2:=#3=#4$}
		\newcommand{\envMThPt}[5][*]{$\ifx d#1 \displaystyle \fi#3#2#4#2#5$}
		\newcommand{\envMThPte}[6][*]{$\ifx d#1 \displaystyle \fi#2#3#4#5#6$}
		\newcommand{\envMF}[5][*]{$\ifx d#1 \displaystyle \fi#2=#3=#4=#5$}
		\newcommand{\envMFPt}[6][*]{$\ifx d#1 \displaystyle \fi#3#2#4#2#5#2#6$}
		
		\newcommand{\envHLineT}[3][*]{ \glcondEnvLineHead{#1} #2&=&#3\glcondEnvLineTail{#1} }
		\newcommand{\envHLineTDef}[3][*]{ \glcondEnvLineHead{#1} #2&:=&#3\glcondEnvLineTail{#1} }

			\newcommand{\envHLine}{\envHLineT}
			\newcommand{\envHLineDef}{\envHLineTDef}

		\newcommand{\envPLine}[2][*]{\glcondEnvLineHead{#1} #2\glcondEnvLineTail{#1}}
		
	\newcommand{\matu}[1][?]{\ifx#1?I\else I_{#1}\fi}
	
	\newcommand{\vlA}[2][n]{\bkAll[#1]{|}{|}{#2}}
				
	\newcommand{\nmSet}[2][n]{\vlA[#1]{#2}}	
		
	\newcommand{\sTx}[2][?]{ \ifx t#1{\tiny #2} \else \ifx s#1{\scriptsize #2} \else \ifx f#1{\footnotesize #2} \else \ifx S#1{\small #2} \else \ifx n#1{\normalsize #2} \else \ifx l#1{\large #2} \else \ifx L#1{\Large #2} \else \ifx R#1{\LARGE #2} \else \ifx h#1{\huge #2} \else \ifx H#1{\Huge #2} \else \ifx ?#1 #2 \else #2 \fi\fi\fi\fi\fi\fi\fi\fi\fi\fi\fi }
	\newcommand{\bfTx}[1]{{\bf#1}}
	
	\newcommand{\ulMTx}[2][?]{\underline{#2}}

	\newcommand{\osMTx}[3][?]{\overset{#3}{#2}}
	\newcommand{\usbMTx}[3][?]{\underset{#2}{\underbrace{#3}}}

		\newcommand{\ulTx}{\ulMTx}

		\newcommand{\osTx}{\osMTx}
		\newcommand{\usbTx}{\usbMTx}
		
	\newcommand{\raTx}[3][?]{\raisebox{#2pt}[0pt][0pt]{\ifx d#1\displaystyle\fi#3}}
	\newcommand{\raMTx}[3][?]{\raisebox{#2pt}[0pt][0pt]{$\ifx d#1\displaystyle\fi#3$}}

	\newcommand{\envCenter}[2][*]{\ifx*#1\begin{center}\else\begin{center}[#1]\fi #2\end{center}}
	\newcommand{\envFlushleft}[2][*]{\ifx*#1\begin{flushleft}\else\begin{flushleft}[#1]\fi #2\end{flushleft}}
	\newcommand{\envFlushright}[2][*]{\ifx*#1\begin{flushright}\else\begin{flushright}[#1]\fi #2\end{flushright}}
		
	\newcommand{\envItemIm}[2][*]{\ifx*#1\begin{itemize}\else\begin{itemize}[#1]\fi #2\end{itemize}}
	\newcommand{\envItemDp}[2][*]{\ifx*#1\begin{description}\else\begin{description}[#1]\fi #2\end{description}}
	\newcommand{\envItemEm}[2][*]{\ifx*#1\begin{enumerate}\else\begin{enfumerate}[#1]\fi #2\end{enumerate}}
		\newcommand{\envItem}{\envItemIm}
		
	\newcommand{\envMultCol}[3][*]{\ifx1#2#3\else\begin{multicols}{#2}\ifx*#1\else\mbox{}\vspace{-#1pt}\fi#3\end{multicols}\fi}
	
\theoremstyle{plain}
\newtheorem{theorem}{Theorem}[section]		
\newtheorem{proposition}[theorem]{Proposition}		
\newtheorem{experiment}[theorem]{Experiment}		
\theoremstyle{definition}
\newtheorem{definition}[theorem]{Definition}		
\newtheorem{problem}[theorem]{Problem}	
\theoremstyle{remark}
\theoremstyle{plain}
\newtheorem{theoremA}{Theorem}[section]
\newtheorem{propositionA}[theoremA]{Proposition}

\theoremstyle{definition}

\newtheorem{algorithmA}[theoremA]{Algorithm}
\newtheorem{processA}[theoremA]{Process}

\theoremstyle{remark}
\newtheorem{remarkA}[theoremA]{Remark}
\theoremstyle{plain}

\theoremstyle{definition}

\theoremstyle{remark}

\allowdisplaybreaks[4]
\numberwithin{equation}{section}

	\newcommand{\lccondBibitem}[3][]{ \if ?#2 \bibitem{#3} \else \bibitem[#2]{#3} \fi}
	\newcommand{\refPaper}[8][?]{
			\lccondBibitem{#1}{#2}
				#3,			
				\emph{#4}, 	
				#5\ 			
				{\bf #6}		
				(#7),			
				#8.			
		}
	
	\newcommand{\refPaperRep}[9][?]{
			\lccondBibitem{#1}{#2}
				#3,			
				\emph{#4}, 	
				#5\ 			
				{\bf #6}		
				(#7),			
				#8			
				; reprinted in #9	
		}
	
	\newcommand{\refBook}[7][?]{
			\lccondBibitem{#1}{#2}
				#3,			
				\emph{#4}, 	
				#5,			
				#6,			
				#7.			
		}
	\newcommand{\refPaperAlm}[5][?]{
			\lccondBibitem{#1}{#2}
				#3,	 		
				\emph{#4}, 	
				#5		
		}

	\newcommand{\etalTx}[2][?]{#2 \emph{et al.}\!}

	\newcommand{\initWd}[2][?]{\ifx q#1`#2'\else\textit{#2}\fi}		
	
	\newcommand{\glcondEnvLineTailPd}[1]{.\ifx*#1 \end{eqnarray*} \else \end{eqnarray} \fi  }
	\newcommand{\glcondEnvLineTailCm}[1]{,\ifx*#1 \end{eqnarray*} \else \end{eqnarray} \fi }
	\newcommand{\prcondEnvEqSpHead}[1]{ \ifx*#1 \begin{equation*}[ERROR] \else \begin{equation}  \label{#1} \fi  }
	\newcommand{\prcondEnvEqSpTail}[1]{\ifx*#1 [ERROR]\end{equation*} \else \end{equation} \fi }

	\newcommand{\envProof}[2][?]{ \par\mbox{}\vspace{-5pt}\\ \ifx?#1\emph{Proof.}\else\emph{Proof of #1.}\fi \ #2 \hfill $\Box$\\ \par}
	
		\newcommand{\envLineTCm}[3][*]{ \glcondEnvLineHead{#1} & &\glcmdEqShift#2\nonumber\\&=&#3 \glcondEnvLineTailCm{#1} }

			\newcommand{\envLineCm}{\envLineTCm}
			
		\newcommand{\envHLineTPd}[3][*]{ \glcondEnvLineHead{#1} #2&=&#3\glcondEnvLineTailPd{#1} }
		\newcommand{\envHLineTDefPd}[3][*]{ \glcondEnvLineHead{#1} #2&:=&#3\glcondEnvLineTailPd{#1} }
		\newcommand{\envHLineTCm}[3][*]{ \glcondEnvLineHead{#1} #2&=&#3\glcondEnvLineTailCm{#1} }
		\newcommand{\envHLineTCmDef}[3][*]{ \glcondEnvLineHead{#1} #2&:=&#3\glcondEnvLineTailCm{#1} }
		\newcommand{\envHLineTCmPt}[4][*]{\glcondEnvLineHead{#1} #3&#2&#4\glcondEnvLineTailCm{#1}}					
		\newcommand{\envHLineTPdPt}[4][*]{\glcondEnvLineHead{#1} #3&#2&#4\glcondEnvLineTailPd{#1}}

			\newcommand{\envHLinePd}{\envHLineTPd}
			\newcommand{\envHLineDefPd}{\envHLineTDefPd}
			\newcommand{\envHLineCm}{\envHLineTCm}
			\newcommand{\envHLineCmDef}{\envHLineTCmDef}
			\newcommand{\envHLineCmPt}{\envHLineTCmPt}
			\newcommand{\envHLinePdPt}{\envHLineTPdPt}

		\newcommand{\envHLineFCmPt}[6][*]{\glcondEnvLineHead{#1} #3&#2&#4\nonumber\\&#2&#5\nonumber \\&#2&#6\glcondEnvLineTailCm{#1}}

		\newcommand{\envHLineCFPd}[5][*]{\glcondEnvLineHead{#1} #2&=&#3,\nonumber\\\glcdH#4&=&#5\glcondEnvLineTailPd{#1}}

		\newcommand{\envHLineCFCmNme}[5][*]{\begin{eqnarray} #2&=&#3,\\\glcdH#4&=&#5 \glcondEnvLineTailCm{?} }
		\newcommand{\envHLineCFNmePd}[5][*]{\begin{eqnarray} #2&=&#3,\\\glcdH#4&=&#5 \glcondEnvLineTailPd{?} }
		\newcommand{\envHLineCFCmDefNme}[5][*]{\begin{eqnarray} #2&:=&#3,\\\glcdH#4&:=&#5 \glcondEnvLineTailCm{?} }
		\newcommand{\envHLineCFDefNmePd}[5][*]{\begin{eqnarray} #2&:=&#3,\\\glcdH#4&:=&#5 \glcondEnvLineTailPd{?} }

		\newcommand{\envHLineCFNmePdPt}[6][*]{\begin{eqnarray}#3&#2&#4,\\\glcdH#5&#2&#6\glcondEnvLineTailPd{?}}
		\newcommand{\envHLineCFCmNmePt}[6][*]{\begin{eqnarray}#3&#2&#4,\\\glcdH#5&#2&#6\glcondEnvLineTailCm{?}}
		\newcommand{\envHLineCFNmePdPte}[7][*]{\begin{eqnarray}#2&#3&#4,\\\glcdH#5&#6&#7\glcondEnvLineTailPd{?}}
		\newcommand{\envHLineCFCmNmePte}[7][*]{\begin{eqnarray}#2&#3&#4,\\\glcdH#5&#6&#7\glcondEnvLineTailCm{?}}

		\newcommand{\envHLineCSNmePd}[7][*]{\begin{eqnarray} #2&=&#3,\\\glcdH#4&=&#5,\\\glcdH#6&=&#7\glcondEnvLineTailPd{?}}
		\newcommand{\envHLineCSDefNmePd}[7][*]{\begin{eqnarray} #2&:=&#3,\\\glcdH#4&:=&#5,\\\glcdH#6&:=&#7\glcondEnvLineTailPd{?}}
		\newcommand{\envHLineCSCmNme}[7][*]{\begin{eqnarray} #2&=&#3,\\\glcdH#4&=&#5,\\\glcdH#6&=&#7\glcondEnvLineTailCm{?}}
		\newcommand{\envHLineCSCmDefNme}[7][*]{\begin{eqnarray} #2&:=&#3,\\\glcdH#4&:=&#5,\\\glcdH#6&:=&#7\glcondEnvLineTailCm{?}}

		\newcommand{\envHLineCSNmePdPt}[8][*]{\begin{eqnarray}#3&#2&#4,\\\glcdH#5&#2&#6,\\\glcdH#7&#2&#8\glcondEnvLineTailPd{?}}
		\newcommand{\envHLineCSCmNmePt}[8][*]{\begin{eqnarray}#3&#2&#4,\\\glcdH#5&#2&#6,\\\glcdH#7&#2&#8\glcondEnvLineTailCm{?}}
		\newcommand{\envHLineCSNmePdPte}[9][*]{\begin{eqnarray}#2&#3&#4,\\\glcdH#5&#6&#7,\\\glcdH#8&#2&#9\glcondEnvLineTailPd{?}}
		\newcommand{\envHLineCSCmNmePte}[9][*]{\begin{eqnarray}#2&#3&#4,\\\glcdH#5&#6&#7,\\\glcdH#8&#2&#9\glcondEnvLineTailCm{?}}
		
		\newcommand{\envHLineCECm}[9][*]{\glcondEnvLineHead{#1} #2&=&#3,\nonumber\\\glcdH#4&=&#5,\nonumber\\\glcdH#6&=&#7,\\\glcdH#8&=&#9\nonumber\glcondEnvLineTailCm{#1}}
		
		\newcommand{\envHLineCENmePd}[9][*]{\begin{eqnarray} #2&=&#3,\\\glcdH#4&=&#5,\\\glcdH#6&=&#7,\\\glcdH#8&=&#9\glcondEnvLineTailPd{?}}
		\newcommand{\envHLineCEDefNmePd}[9][*]{\begin{eqnarray} #2&:=&#3,\\\glcdH#4&:=&#5,\\\glcdH#6&:=&#7,\\\glcdH#8&:=&#9\glcondEnvLineTailPd{?}}
		\newcommand{\envHLineCECmNme}[9][*]{\begin{eqnarray} #2&=&#3,\\\glcdH#4&=&#5,\\\glcdH#6&=&#7,\\\glcdH#8&=&#9\glcondEnvLineTailCm{?}}
		\newcommand{\envHLineCECmDefNme}[9][*]{\begin{eqnarray} #2&:=&#3,\\\glcdH#4&:=&#5,\\\glcdH#6&:=&#7,\\\glcdH#8&:=&#9\glcondEnvLineTailCm{?}}

			\newcommand{\pccondPaForPar}[1]{\ifx#1p \\\glcdH \fi}
			\newcommand{\pccondPaForNonnum}[1]{\ifx#1* \nonumber \fi}
			\newcommand{\HLineCTCm}[3][?]{#2&=&#3, \nonumber\pccondPaForPar{#1}}

			\newcommand{\HLineCECm}[9][?]{#2&=&#3,\nonumber\\\glcdH#4&=&#5,\nonumber\\\glcdH#6&=&#7,\nonumber\\\glcdH#8&=&#9,\pccondPaForPar{#1}}
			\newcommand{\HLineCEPd}[9][?]{#2&=&#3,\nonumber\\\glcdH#4&=&#5,\nonumber\\\glcdH#6&=&#7,\nonumber\\\glcdH#8&=&#9.\pccondPaForPar{#1}}

		\newcommand{\envPLinePd}[2][*]{\glcondEnvLineHead{#1} #2\glcondEnvLineTailPd{#1}}
		\newcommand{\envPLineCm}[2][*]{\glcondEnvLineHead{#1} #2\glcondEnvLineTailCm{#1}}
	
		\newcommand{\envOTLinePd}[4][*]{\glcondEnvLineHead{#1} #2\lnAP{=}#3\lnP{=}#4.\glcondEnvLineTail{#1}}

		\newcommand{\envOTLinePdPt}[5][*]{\glcondEnvLineHead{#1} #3\lnAP{#2}#4\lnP{#2}#5\glcondEnvLineTailPd{#1}}

			\newcommand{\envOTLineThPdPt}{\envOTLinePdPt}

		\newcommand{\envOFLineCm}[5][*]{\glcondEnvLineHead{#1} #2\lnAP{=}#3\lnP{=}#4\lnP{=}#5,\glcondEnvLineTail{#1}}

		\newcommand{\envMOCm}[2][*]{$\ifx d#1 \displaystyle \fi#2$,}
		\newcommand{\envMOPd}[2][*]{$\ifx d#1 \displaystyle \fi#2$.}
		
		\newcommand{\envMTCm}[3][*]{$\ifx d#1 \displaystyle \fi#2=#3$,}
		\newcommand{\envMTPd}[3][*]{$\ifx d#1 \displaystyle \fi#2=#3$.}
		\newcommand{\envMTCmDef}[3][*]{$\ifx d#1 \displaystyle \fi#2:=#3$,}
		\newcommand{\envMTDefPd}[3][*]{$\ifx d#1 \displaystyle \fi#2:=#3$.}
		\newcommand{\envMTCmPt}[4][*]{$\ifx d#1 \displaystyle \fi#3#2#4$,}
		\newcommand{\envMTPdPt}[4][*]{$\ifx d#1 \displaystyle \fi#3#2#4$.}
			\newcommand{\envMCm}{\envMTCm}
			\newcommand{\envMPd}{\envMTPd}

			\newcommand{\envMPdPt}{\envMTPdPt}
		\newcommand{\envMThCm}[4][*]{$\ifx d#1 \displaystyle \fi#2=#3=#4$,}
		\newcommand{\envMThPd}[4][*]{$\ifx d#1 \displaystyle \fi#2=#3=#4$.}
		\newcommand{\envMThCmPt}[5][*]{$\ifx d#1 \displaystyle \fi#3#2#4#2#5$,}
		\newcommand{\envMThPdPt}[5][*]{$\ifx d#1 \displaystyle \fi#3#2#4#2#5$.}
		\newcommand{\envMFCm}[5][*]{$\ifx d#1 \displaystyle \fi#2=#3=#4=#5$,}
		\newcommand{\envMFPd}[5][*]{$\ifx d#1 \displaystyle \fi#2=#3=#4=#5$.}
		\newcommand{\envMFCmPt}[6][*]{$\ifx d#1 \displaystyle \fi#3#2#4#2#5#2#6$,}
		\newcommand{\envMFPdPt}[6][*]{$\ifx d#1 \displaystyle \fi#3#2#4#2#5#2#6$.}
		\newcommand{\envMFiCm}[6][*]{$\ifx d#1 \displaystyle \fi#2=#3=#4=#5=#6$,}
		\newcommand{\envMFiPd}[6][*]{$\ifx d#1 \displaystyle \fi#2=#3=#4=#5=#6$.}
		\newcommand{\envMFiCmPt}[7][*]{$\ifx d#1 \displaystyle \fi#3#2#4#2#5#2#6#2#7$,}
		\newcommand{\envMFiPdPt}[7][*]{$\ifx d#1 \displaystyle \fi#3#2#4#2#5#2#6#2#7$.}

		\newcommand{\prcondHLCPNm}{\hspace{-1pt}}
		\newcommand{\envHLineCFCmNm}[5][*]{ \begin{equation}\begin{split} \ifx*#1 \text{[ERROR;need label name]} \else \label{#1} \fi #2&\prcondHLCPNm\lnP{=}\prcondHLCPNm#3,\\#4&\prcondHLCPNm\lnP{=}\prcondHLCPNm#5, \end{split}\end{equation} }
		\newcommand{\envHLineCFNm}[5][*]{ \begin{equation}\begin{split} \ifx*#1 \text{[ERROR;need label name]} \else \label{#1} \fi #2&\prcondHLCPNm\lnP{=}\prcondHLCPNm#3\\#4&\prcondHLCPNm\lnP{=}\prcondHLCPNm#5, \end{split}\end{equation} }
		\newcommand{\envHLineCFNmPd}[5][*]{ \begin{equation}\begin{split} \ifx*#1 \text{[ERROR;need label name]} \else \label{#1} \fi #2&\prcondHLCPNm\lnP{=}\prcondHLCPNm#3,\\#4&\prcondHLCPNm\lnP{=}\prcondHLCPNm#5. \end{split}\end{equation} }
		\newcommand{\envHLineCFCmDefNm}[5][*]{ \begin{equation}\begin{split} \ifx*#1 \text{[ERROR;need label name]} \else \label{#1} \fi #2&\prcondHLCPNm\lnP{:=}\prcondHLCPNm#3,\\#4&\prcondHLCPNm\lnP{:=}\prcondHLCPNm#5, \end{split}\end{equation} }
		\newcommand{\envHLineCFDefNm}[5][*]{ \begin{equation}\begin{split} \ifx*#1 \text{[ERROR;need label name]} \else \label{#1} \fi #2&\prcondHLCPNm\lnP{:=}\prcondHLCPNm#3\\#4&\prcondHLCPNm\lnP{:=}\prcondHLCPNm#5, \end{split}\end{equation} }
		\newcommand{\envHLineCFDefNmPd}[5][*]{ \begin{equation}\begin{split} \ifx*#1 \text{[ERROR;need label name]} \else \label{#1} \fi #2&\prcondHLCPNm\lnP{:=}\prcondHLCPNm#3,\\#4&\prcondHLCPNm\lnP{:=}\prcondHLCPNm#5. \end{split}\end{equation} }
		\newcommand{\envHLineCSCmNm}[7][*]{ \begin{equation}\begin{split} \ifx*#1 \text{[ERROR;need label name]} \else \label{#1} \fi #2&\prcondHLCPNm\lnP{=}\prcondHLCPNm#3,\\#4&\prcondHLCPNm\lnP{=}\prcondHLCPNm#5,\\#6&\prcondHLCPNm\lnP{=}\prcondHLCPNm#7 \end{split}\end{equation} }
		\newcommand{\envHLineCSNm}[7][*]{ \begin{equation}\begin{split} \ifx*#1 \text{[ERROR;need label name]} \else \label{#1} \fi #2&\prcondHLCPNm\lnP{=}\prcondHLCPNm#3\\#4&\prcondHLCPNm\lnP{=}\prcondHLCPNm#5\\#6&\prcondHLCPNm\lnP{=}\prcondHLCPNm#7 \end{split}\end{equation} }
		\newcommand{\envHLineCSNmPd}[7][*]{ \begin{equation}\begin{split} \ifx*#1 \text{[ERROR;need label name]} \else \label{#1} \fi #2&\prcondHLCPNm\lnP{=}\prcondHLCPNm#3,\\#4&\prcondHLCPNm\lnP{=}\prcondHLCPNm#5,\\#6&\prcondHLCPNm\lnP{=}\prcondHLCPNm#7. \end{split}\end{equation} }
		\newcommand{\envHLineCSCmDefNm}[7][*]{ \begin{equation}\begin{split} \ifx*#1 \text{[ERROR;need label name]} \else \label{#1} \fi #2&\prcondHLCPNm\lnP{:=}\prcondHLCPNm#3,\\#4&\prcondHLCPNm\lnP{:=}\prcondHLCPNm#5,\\#6&\prcondHLCPNm\lnP{:=}\prcondHLCPNm#7 \end{split}\end{equation} }
		\newcommand{\envHLineCSDefNm}[7][*]{ \begin{equation}\begin{split} \ifx*#1 \text{[ERROR;need label name]} \else \label{#1} \fi #2&\prcondHLCPNm\lnP{:=}\prcondHLCPNm#3\\#4&\prcondHLCPNm\lnP{:=}\prcondHLCPNm#5\\#6&\prcondHLCPNm\lnP{:=}\prcondHLCPNm#7 \end{split}\end{equation} }
		\newcommand{\envHLineCSDefNmPd}[7][*]{ \begin{equation}\begin{split} \ifx*#1 \text{[ERROR;need label name]} \else \label{#1} \fi #2&\prcondHLCPNm\lnP{:=}\prcondHLCPNm#3,\\#4&\prcondHLCPNm\lnP{:=}\prcondHLCPNm#5,\\#6&\prcondHLCPNm\lnP{:=}\prcondHLCPNm#7. \end{split}\end{equation} }
		\newcommand{\envHLineCECmNm}[9][*]{ \begin{equation}\begin{split} \ifx*#1 \text{[ERROR;need label name]} \else \label{#1} \fi #2&\prcondHLCPNm\lnP{=}\prcondHLCPNm#3,\\#4&\prcondHLCPNm\lnP{=}\prcondHLCPNm#5,\\#6&\prcondHLCPNm\lnP{=}\prcondHLCPNm#7,\\#8&\prcondHLCPNm\lnP{=}\prcondHLCPNm#9,  \end{split}\end{equation} }
		\newcommand{\envHLineCENm}[9][*]{ \begin{equation}\begin{split} \ifx*#1 \text{[ERROR;need label name]} \else \label{#1} \fi #2&\prcondHLCPNm\lnP{=}\prcondHLCPNm#3\\#4&\prcondHLCPNm\lnP{=}\prcondHLCPNm#5\\#6&\prcondHLCPNm\lnP{=}\prcondHLCPNm#7\\#8&\prcondHLCPNm\lnP{=}\prcondHLCPNm#9  \end{split}\end{equation} }
		\newcommand{\envHLineCENmPd}[9][*]{ \begin{equation}\begin{split} \ifx*#1 \text{[ERROR;need label name]} \else \label{#1} \fi #2&\prcondHLCPNm\lnP{=}\prcondHLCPNm#3,\\#4&\prcondHLCPNm\lnP{=}\prcondHLCPNm#5,\\#6&\prcondHLCPNm\lnP{=}\prcondHLCPNm#7,\\#8&\prcondHLCPNm\lnP{=}\prcondHLCPNm#9.  \end{split}\end{equation} }
		\newcommand{\envHLineCECmDefNm}[9][*]{ \begin{equation}\begin{split} \ifx*#1 \text{[ERROR;need label name]} \else \label{#1} \fi #2&\prcondHLCPNm\lnP{:=}\prcondHLCPNm#3,\\#4&\prcondHLCPNm\lnP{:=}\prcondHLCPNm#5,\\#6&\prcondHLCPNm\lnP{:=}\prcondHLCPNm#7,\\#8&\prcondHLCPNm\lnP{:=}\prcondHLCPNm#9,  \end{split}\end{equation} }
		\newcommand{\envHLineCEDefNm}[9][*]{ \begin{equation}\begin{split} \ifx*#1 \text{[ERROR;need label name]} \else \label{#1} \fi #2&\prcondHLCPNm\lnP{:=}\prcondHLCPNm#3\\#4&\prcondHLCPNm\lnP{:=}\prcondHLCPNm#5\\#6&\prcondHLCPNm\lnP{:=}\prcondHLCPNm#7\\#8&\prcondHLCPNm\lnP{:=}\prcondHLCPNm#9  \end{split}\end{equation} }
		\newcommand{\envHLineCEDefNmPd}[9][*]{ \begin{equation}\begin{split} \ifx*#1 \text{[ERROR;need label name]} \else \label{#1} \fi #2&\prcondHLCPNm\lnP{:=}\prcondHLCPNm#3,\\#4&\prcondHLCPNm\lnP{:=}\prcondHLCPNm#5,\\#6&\prcondHLCPNm\lnP{:=}\prcondHLCPNm#7,\\#8&\prcondHLCPNm\lnP{:=}\prcondHLCPNm#9.  \end{split}\end{equation} }
	\newcommand{\envMLineTPd}[3][*]{ \ifx*#1 \begin{multline*} #2\lnP{=}#3.\end{multline*} \else \begin{multline} \label{#1} #2\lnP{=}#3.\end{multline} \fi }
	\newcommand{\envMLineTCm}[3][*]{ \ifx*#1 \begin{multline*} #2\lnP{=}#3,\end{multline*} \else \begin{multline} \label{#1} #2\lnP{=}#3,\end{multline} \fi }
	\newcommand{\envMLineTDefPd}[3][*]{ \ifx*#1 \begin{multline*} #2\lnP{:=}#3.\end{multline*} \else \begin{multline} \label{#1} #2\lnP{:=}#3.\end{multline} \fi }
	\newcommand{\envMLineTCmDef}[3][*]{ \ifx*#1 \begin{multline*} #2\lnP{:=}#3,\end{multline*} \else \begin{multline} \label{#1} #2\lnP{:=}#3,\end{multline} \fi }

	\DeclareFontEncoding{OT2}{}{}
	\DeclareFontSubstitution{OT2}{cmr}{m}{n}
	\DeclareFontFamily{OT2}{cmr}{\hyphenchar\font45}
	\DeclareFontShape{OT2}{cmr}{m}{n}{<5><6><7><8><9>gen*wncyr <10><10.95><12><14.4><17.28><20.74><24.88>wncyr10}{}
	\DeclareFontShape{OT2}{cmr}{b}{n}{<5><6><7><8><9>gen*wncyb<10><10.95><12><14.4><17.28><20.74><24.88>wncyb10}{}
	\DeclareMathAlphabet{\mathcyr}{OT2}{cmr}{m}{n}
	\DeclareMathAlphabet{\mathcyb}{OT2}{cmr}{b}{n}
	\SetMathAlphabet{\mathcyr}{bold}{OT2}{cmr}{b}{n}

	\newcommand{\wgt}[2][n]{\Fc[#1]{{\mathrm{w}}}{#2}}	
	\newcommand{\dpt}[2][n]{\Fc[#1]{{\mathrm{d}}}{#2}}

	\newcommand{\shh}{*}
	\newcommand{\shs}{\mathcyr{sh}}	
	\newcommand{\shH}[1][\,]{#1\shh#1}
	\newcommand{\shS}[1][\;]{#1\shs#1}
	
	\newcommand{\spZN}{\mathcal{Z}}
	\newcommand{\spZrN}{\rLt{\spZN}}		
	\newcommand{\spZO}[2][?]{\if ?#1 \spZN_{#2} \else \spZN_{#2}^{#1} \fi}
	\newcommand{\spZrO}[2][?]{\if ?#1 \spZrN_{#2} \else \spZrN_{#2}^{#1} \fi}		
	\newcommand{\spZT}[3][?]{\if ?#1 \spZN_{#2,#3} \else \spZN_{#2,#3}^{#1} \fi}	
	\newcommand{\spZrT}[3][?]{\if ?#1 \spZrN_{#2,#3} \else \spZrN_{#2,#3}^{#1} \fi}	
	\newcommand{\spHN}[1][?]{\if ?#1 \mathfrak{H} \else \mathfrak{H}_{#1} \fi}	
	\newcommand{\spHrN}[1][?]{\if ?#1 \spHN^1 \else \spHN^1_{#1} \fi} 		
	\newcommand{\spHrrN}[1][?]{\if ?#1 \spHN^0 \else \spHN^0_{#1} \fi}		

	\newcommand{\evZO}[2][n]{\Fc[#1]{Z}{#2}}

	\newcommand{\cdmZO}[2][?]{d_{#2}}		

	\newcommand{\ztN}[1][?]{\zeta}				\newcommand{\ztO}[2][n]{\Fc[#1]{\zeta}{#2}}

	\newcommand{\GztN}[1][?]{\bLt{\zeta}}

		\newcommand{\zt}{\ztO}				

\input{letters.tex}

	\newcommand{\FLEm}[1][?]{\mathtt{#1\_M.py}}
	\newcommand{\FLEc}[1][?]{\mathtt{#1\_C.py}}

	\newcommand{\spZqT}{\spZrT}

	\newcommand{\cdmBZqT}[3][?]{\letD^{\letB}_{#2,#3}}		
	\newcommand{\cdmBZMO}[2][?]{\letD^{\letB,\letMJPO}_{#2}}		
	\newcommand{\cdmBZKO}[2][?]{\letD^{\letB,\letKNT}_{#2}}		
	\newcommand{\cdmBZBO}[2][?]{\letD^{\letB,\bullet}_{#2}}		
	\newcommand{\spEO}[2][?]{\mathcal{E}_{#2}}			
	\newcommand{\bspHN}[1][?]{\spHN^{\setF[2]}}		
	\newcommand{\bspUN}[1][?]{\mathcal{H}^{\letB}}		
	\newcommand{\bspUO}[2][?]{\bspUN_{#2}}
	\newcommand{\bspUT}[3][?]{\bspUN_{#2,#3}}	
	\newcommand{\bspEO}[2][?]{\mathcal{E}^{\letB}_{#2}}			
	\newcommand{\bspEBO}[2][?]{\mathcal{E}^{\letB,\bullet}_{#2}}			
	\newcommand{\bspZN}[1][?]{\spZN^{\letB}}	
	\newcommand{\bspZO}[2][?]{\spZO[\letB]{#2}}	
	\newcommand{\bspZT}[3][?]{\spZT[\letB]{#2}{#3}}	%
	\newcommand{\bspZqT}[3][?]{\spZrT[\letB]{#2}{#3}}	
	\newcommand{\bspEqT}[3][?]{\rLt{\mathcal{E}}^{\letB}_{#2,#3}}			
	\newcommand{\ispHN}[1][?]{\spHN^{\setZ}}		
	\newcommand{\ispHrrN}[1][?]{\spHN^{\setZ,0}}	
	\newcommand{\ispEO}[2][?]{\mathcal{E}^{\setZ}_{#2}}		

	\newcommand{\bztO}[2][n]{\pwFc[#1]{\zeta}{\letB}{#2}}	
	\newcommand{\bztqO}[2][n]{\pwFc[#1]{\rLt{\zeta}}{\letB}{#2}}	
	\newcommand{\bztSO}[2][n]{\ipFc[#1]{\zeta}{\shs}{\letB}{#2}}	
	\newcommand{\bZtO}[2][n]{\pwFc[#1]{\eta}{\letB}{#2}}	
		\newcommand{\bzt}{\bztO}
		\newcommand{\bztS}{\bztSO}
		\newcommand{\bZt}{\bZtO}
		
	\newcommand{\bevHO}[2][n]{\pwFc[#1]{H}{\letB}{#2}}		%
	\newcommand{\bevZO}[2][n]{\pwFc[#1]{\evZO[]{}}{\letB}{#2}}		%
	\newcommand{\rlDO}[2][n]{\Fc[#1]{\mathsf{ds}}{#2}}		

	\newcommand{\hmRLett}{\mathrm{reg}}
	\newcommand{\hmsRO}[2][n]{\idFc[#1]{\hmRLett}{\shs}{#2}}	
	\newcommand{\hmhRO}[2][n]{\idFc[#1]{\hmRLett}{\shh}{#2}}
	
	\newcommand{\lnCLett}{\mathrm{can}}
	\newcommand{\lnCO}[2][n]{\pwFc[#1]{\lnCLett}{\letB}{#2}}	
	
	\newcommand{\midE}{\varnothing}		
	\newcommand{\midO}[1][?]{\ifx #1?\bfTx{1}\else\bfTx{1}_{#1}\fi}		
	\newcommand{\midsN}[1][?]{\bfTx{I}}	
	\newcommand{\midsO}[2][?]{\midsN_{#2}}	
	\newcommand{\midsT}[3][?]{\midsN_{#2,#3}}	

	\newcommand{\midsHN}[1][?]{\midsN^{H}}	
	\newcommand{\midsHO}[2][?]{\midsHN_{#2}}		
	\newcommand{\midsHT}[3][?]{\midsHN_{#2,#3}}		
	\newcommand{\emids}[1][?]{\eLt{\midsN}}	
	\newcommand{\emidsN}[1][?]{\bfTx{\eLt{I}}}	
	\newcommand{\emidsO}[2][?]{\emidsN_{#2}}	
	\newcommand{\pidsN}[1][?]{\bfTx{PI}}	
	\newcommand{\epidsN}[1][?]{\eLt{\pidsN}}	
	\newcommand{\pidsO}[2][?]{\pidsN_{#2}}
	\newcommand{\epidsO}[2][?]{\epidsN_{k}}
	\newcommand{\epidsEO}[2][?]{\epidsN_{k}^{\mathrm{EDS}}}	
	\newcommand{\epidsMO}[2][?]{\epidsN_{k}^{\letMJPO}}	
	\newcommand{\epidsKO}[2][?]{\epidsN_{k}^{\letKNT}}	
	\newcommand{\epidsBO}[2][?]{\epidsN_{k}^{\bullet}}	

	\newcommand{\letSPD}{\mathtt{Sh}}	
	\newcommand{\letHPD}{\mathtt{St}}
	\newcommand{\spdO}[2][n]{\Fc[#1]{\letSPD}{#2}}
	\newcommand{\hpdO}[2][n]{\Fc[#1]{\letHPD}{#2}}
	\newcommand{\spdRO}[2][n]{\idFc[#1]{\letSPD}{\shs}{#2}}		
	\newcommand{\hpdRO}[2][n]{\idFc[#1]{\letHPD}{\shs}{#2}}		
	
	\newcommand{\soSD}[1][1]{\bigtriangleup}	

	\newcommand{\sysN}[1][?]{\mathrm{S}}

	\newcommand{\sbsFO}[2][?]{\idFc{s}{\mathrm{min}}{#2}}		
	\newcommand{\sbssPN}[1][?]{P}		
	\newcommand{\sbssDN}[1][?]{D}		

	\newcommand{\sbssDT}[3][?]{\sbssDN_{#2,#3}}

	\newcommand{\coeFO}[2][?]{\idFc{c}{\mathrm{min}}{#2}}		

	\newcommand{\cmbsN}[1][*]{\mathcal{K}^{#1}}
	\newcommand{\cmbsO}[2][*]{\cmbsN[#1]_{#2}}	

	\newcommand{\cndBO}[2][?]{\ifx1#2(a)\else\ifx2#2(b)\else!ERROR!\fi\fi}		
	\newcommand{\cndEO}[2][?]{(E#2)}
	\newcommand{\cndEdO}[2][?]{(e#2)}
	
	\newcommand{\ltm}[2][n]{\Fc[#1]{L}{#2}}

\allowdisplaybreaks[4]
	\title{\mainTitle}
	\author{
			\authorName\thanks{\organizationNameFst, \placeAddressFst}
		}
	
	\date{}

\begin{document}
\maketitle
\renewcommand{\thefootnote}{\fnsymbol{footnote}}
\footnote[0]{e-mail : \emailAddressFst}
\footnote[0]{MSC-class: \MSCname}
\renewcommand{\thefootnote}{\arabic{footnote}}\setcounter{footnote}{0}
\vPack[30]

\begin{abstract}
The formal multiple zeta space we consider with a computer 
	is 
	an $\setF[2]$-vector space 
	generated by $2^{k-2}$ formal symbols for a given weight $k$,
	where 
	the symbols satisfy binary extended double shuffle relations.			
Up to weight $k=22$,
	we compute  
	the dimensions of the formal multiple zeta spaces, 
	and 
	verify 
	the dimension conjecture 
	on original extended double shuffle relations of real multiple zeta values.	
Our computations   
	adopt 
	Gaussian forward elimination
	and
	give 	
	information for spaces filtered by depth.
We can observe that 
	the dimensions of the depth-graded formal multiple zeta spaces 
	have 
	a Pascal triangle pattern 
	expected by the Hoffman mult-indices.
\end{abstract}

\section{Introduction} \label{sectOne}
The space generated by multiple zeta values (MZVs for short)	
	has been elucidated theoretically and numerically in recent years,
	but
	its structure remains mysterious.
In this paper,
	we shed light on  
	a formal space generated by binary analogs of MZVs 
	by computer experiments
	for unraveling both of the original and formal spaces.	
	
Let $\setN$
	denote the set of positive integers.
The MZV is a real number
	that belongs to  
	an image of a function (customarily denoted by $\zt[]{}$)
	whose domain 
	is 
	\envHLineCm[1_PL_DefMIds]
	{
		\midsN
	}
	{	
		\UnO{r \geq 0} \SetT{ \bfK_r=(k_1,k_2,\ldots,k_r) \in \setN^r}{ k_1 \geq 2 }	
	}
	where 
	$\bfK_0=\midE$ is the empty \initWd{mult-index}
	and
	$\zt{\midE}=1$.
We call $\wgt{\bfK_r}=k_1+\cdots+k_r$ and $\dpt{\bfK_r}=r$ the weight and depth,	
	respectively.
The function $\zt[]{}$ 
	has two definitions 
	by the iterated integral and nested summation,
	which endow the $\setQ$-vector space $\spZN$ spanned by MZVs with abundant linear relations.
Euler \cite{Euler1776},
	who solved the Basel problem $\zt{2}=\frc[s]{6}{\pi^2}$ and advanced the case $r=1$,
	also studied the case $r=2$.	

Zagier \cite{Zagier94} conjectured\footnote{
Zagier noted the conjectures were made after many discussions with Drinfel'd, Kontsevich and Goncharov.
}
	that
	$\spZN$ is graded by weight
	and
	the dimensions of graded pieces are expressed in terms of a Fibonacci-like sequence. 
Let $\midsO{k}$ 
	be 
	the subset consisting of mult-indices of weight $k$, 		
	and
	let
	$\spZO{k}$ be the subspace 
	spanned by MZVs 
	in $\zt{\midsO{k}} = \SetT{\zt{\bfK}}{ \bfK \in \midsO{k}}$.	
The dimension conjecture is 
	\envHLineCmPt[1_PL_CnjDimW]{\osTx{=}{?}}
	{
		\dimT{\setQ}{\spZO{k}}
	}{
		\cdmZO{k}
	}
	where 
	$\cdmZO{k}=\cdmZO{k-2}+\cdmZO{k-3}$ $(k\geq3)$,
	$\cdmZO{0}=\cdmZO{2}=1$
	and
	$\cdmZO{1}=0$.
These integers fit together into the generating series 	
	\envHLinePd[1_PL_EqCDmZ]
	{
		\Sm{k \geq 0} \cdmZO{k} X^k
	}
	{
		\frc{1-(X^2+X^3)}{1}
	}
	
The ultimate upper bound theorem 
	(i.e., \envMPt{\leq}{\dimT{\setQ}{\spZO{k}}}{\cdmZO{k}})
	was established independently by Goncharov \cite{DG05,Goncharov02} and Terasoma \cite{Terasoma02}.
Brown \cite{Brown12} furthermore proved that
	$\spZO{k}$ 
	is generated by 	
	MZVs in $\zt{\midsHO{k}}$,
	where
	$\midsHO{k}$ is the set of Hoffman mult-indices of weight $k$:
	\envHLinePd[1_PL_DefHMIdsW]
	{
		\midsHO{k}
	}
	{
		\SetT{ \bfK=(k_1,\ldots,k_r) \in \midsO{k}}{ k_i \in \SetO{2,3}}
	}
Hoffman \cite{Hoffman97} conjectured 
	$\zt{\midsHO{k}}$ is a basis of $\spZO{k}$,
	which 
	would imply 	
	the dimension conjecture
	because 	
	the same recurrence relation 
	$\nmSet{\midsHO{k}}=\nmSet{\midsHO{k-2}}+\nmSet{\midsHO{k-3}}$ 
	holds by a simple count of the number of $2$'s and $3$'s.
Umezawa \cite{Umezawa19} also 
	suggested a basis conjecture 
	in terms of iterated log-sine integrals,
	in which 
	sets of mult-indices different from $\midsHO{k}$ are used.
Because of the difficulty to show the independence between MZVs,
	no non-trivial lower bounds are known.

By the upper bound theorem,
	it is natural to ask that
	what sorts of relations are needed to reduce the number of generators of $\spZO{k}$ to $\cdmZO{k}$.
There are 	several conjectural candidates: 
	e.g., \cite{Drinfeld91,Furusho03,HS19,KY18,Kawashima09}.	
In particular,
	the extended double shuffle (EDS) relations \cite{IKZ06,Racinet02} known from early on 
	are often selected 
	for experimentally attacking this question,	
	because
	they are easier to write down		
	and
	included in the other candidates except Kawasima's \cite{Kawashima09}.		
Minh and Petitot \cite{MP00} 	
	verified 
	that
	the class of EDS relations is a right candidate up to weight $10$,
	\etalTx{Bigotte} \cite{BJOP02} verified it up to weight $12$, 
	\etalTx{Minh} \cite{MJPO00} verified it up to weight $16$,\!\footnote{
This experimental result was announced in their private communication (see \cite[Section\,1]{KNT08}).
}	
	\etalTx{Espie} \cite{ENR02} verified it up to weight $19$, 	
	and
	\etalTx{Kaneko} \cite{KNT08} verified it up to weight $20$ that seems to be the latest record.
The first two experiments are by the Gr\"obner basis method,
	and
	the last three ones are by the vector space (or matrix) method.
The fourth one of \cite{ENR02}
	was executed 
	under	
	modulo rational multiples of powers of $\zt{2}$,
	or module $\setQ[\zt{2}]$.		

The first purpose of this paper is to improve the record to weight $k=22$.
For this,
	we consider 
	an $\setF[2]$-vector space $\bspZO{k}$
	instead of the $\setQ$-vector space $\spZO{k}$:
	roughly speaking,
	$\bspZO{k}$
	is
	generated 
	by \initWd{binary multiple zeta symbols} $\bztO{\bfK}$ $(\bfK \in \midsO{k})$ (binary MZSs for short),
	where 
	$\bztO{\bfK}$ satisfy 
	\initWd{binary EDS relations}
	that are 
	obtained from original EDS relations 
	after 
	the modulo $2$ arithmetic to integer coefficients.  
(Exact definitions of the binary analogs 
	in this section
	will be stated in the next section.)	
We will verify 	
	$\bztO{\midsHO{k}}$ is a basis of $\bspZO{k}$	
	and
	\envMPd	
	{
		\dimT{\setF[2]}{\bspZO{k}}
	}{ 
		\cdmZO{k}
	}
Our calculation results break 
	the record
	because 
	\envMPt{\leq}{ 
		\dimT{\setQ}{\spZO{k}}
	}{
		\dimT{\setF[2]}{\bspZO{k}}
	}	
	(as will be mentioned in \refSect{sectThree}).
The space $\bspZO{k}$
	reduces 
	the computation cost		
	since 	%
	$\setF[2]$ is the binary and simplest finite field.		
The field $\setF[2]$ makes it easy 
	to apply useful techniques in computer 	
	since 
	$\setF[2]$ is compatible with the Boolean datatype:		
	in fact,
	we will employ a conflict based algorithm discussed in \cite{MS18}		
	for a fast Gaussian forward elimination.

The second and main purpose is 
	to observe  
	a Pascal triangle pattern in $\bspZO{k}$ 
	from the viewpoint of a direct sum decomposition,
	\envHLineCmPt[1_PL_CgDSDZ]{\cong}
	{
		\bspZO{k}  
	}
	{
		\bspZqT{k}{k-1} \OPlus \cdots \OPlus  \bspZqT{k}{0}		
	}	
	where
	$\bspZqT{k}{r}$ 
	are
	quotient spaces
	defined 
	by means of depth filtration:	%
	the descending chain
	$\bspZT{k}{k-1} \supset \cdots \supset \bspZT{k}{0}$ 
	is used for $\bspZqT{k}{r} = \bspZT{k}{r} / \bspZT{k}{r-1}$,
	where
	$\bspZT{k}{r}$ 
	are 	
	the subspaces spanned by 
	binary MZSs of weight $k$ and depth at most $r$.	
We define  
	$\midsT{k}{r} = \SetT{ \bfK \in \midsO{k}}{ \dpt{\bfK}=r}$	
	and		
	\envHLineCm[1_PL_DefHMIdsWD]
	{
		\midsHT{k}{r}
	}
	{
		\midsHO{k} \cap \midsT{k}{r}		
	}	
	with 	
	\envMPd{ 
		\cdmBZqT{k}{r}		
	}{
		\nmSet{\midsHT{k}{r}}
	}
We denote by $\bztqO{\bfK}$ 
	the canonical image\footnote{
	We use the same notation $\bztqO[]{}$ 
	for all canonical images in the quotient spaces $\bspZqT{k}{r}$ $(k > r \geq 0)$.	
	There should be no confusion 
	because the quotient space under consideration is clear from context.
}
	of $\bztO{\bfK}$ in $\bspZqT{k}{r}$ for any $\bfK \in \midsT{k}{r}$.
Up to weight $k=22$,	
	we will verify
	$\bztqO{\midsHT{k}{r}}$ is a basis of $\bspZqT{k}{r}$		
	and
	\envMPd		
	{
		\dimT{\setF[2]}{\bspZqT{k}{r}}
	}{ 
		\cdmBZqT{k}{r}		
	}	
Counting the number of $2$'s and $3$'s
	implies 
	that
	the double sequence 
	$\bkR{\cdmBZqT{k}{r}}$ 
	satisfies a recurrence relation 
	with a Pascal triangle pattern:	 
	$\cdmBZqT{k}{r} = \cdmBZqT{k-2}{r-1} + \cdmBZqT{k-3}{r-1}$ $(k\geq3, r \geq 1)$,
	$\cdmBZqT{0}{0}=\cdmBZqT{2}{1}=1$
	and
	$\cdmBZqT{k}{r}=0$ for other $k$ and $r$,
	or equivalently,
	\envHLinePd[1_PL_EqCDmBZq]
	{
		\Sm{k,r \geq 0} \cdmBZqT{k}{r} X^k Y^r
	}
	{
		\frc{1-(X^2+X^3)Y}{1}
	}
More precisely,
	$\cdmBZqT{k}{r} = \bnm{k-2r}{r}$
	since  
	the integers $P_{r,k}=\cdmBZqT{k+2r}{r}$
	satisfy the same recurrence relation as the binomial coefficients $\bnm{k}{r}$.
As expected from \refEq{1_PL_CgDSDZ},
	the formula \refEq{1_PL_EqCDmBZq}
	specializes to \refEq{1_PL_EqCDmZ} 
	upon $Y=1$. 

We also try 
	experiments 
	on parts of EDS relations,
	\initWd[q]{$\letKNT$} and \initWd[q]{$\letMJPO$} relations, 	
	which 
	are expected to be alternatives to EDS
	and
	actually employed in \cite{KNT08,MJPO00} for verification,
	respectively.
Unlike the case in $\spZO{k}$,	
	those relations do not suffice to give all relations in $\bspZO{k}$,
	but 
	we can find 
	a quasi Fibonacci-like rule 
	in dimensions of spaces defined by $\letMJPO$ relations.

The idea of the depth filtration in \refEq{1_PL_CgDSDZ}		
	was conceived	
	by Broadhurst and Kreimer \cite{BK97}
	to
	propose a refinement of the dimension conjecture. 
Their conjecture indicates 
	two interesting facts 	
	in the $\setQ$-vector spaces of MZVs 
	graded by both weight and depth:	
	{(i)}
	modular forms influence the structure through 
	quotient spaces $\spZqT{k}{r}$
	defined by
	the $\setQ$-version of \refEq{1_PL_CgDSDZ};
	and 
	{(ii)}
	the Hoffman values $\zt{\bfK}$ $(\bfK \in \midsHO{k})$ are irrelevant to the structure	
	in the sense that
	most of the values vanish in the graded pieces of same depth.	
In terms of the generating series,	
	the conjecture is 
	\envHLineCmPt[1_PL_EqCDmZqBK]{\osTx{=}{?}}
	{
		\Sm{k,r \geq 0} \dimO{\spZqT{k}{r}} X^k Y^r
	}
	{
		\frc{1 - \Fc{O}{X} Y + \Fc{S}{X} Y^2 (1 -Y^2)}{ 1 + \Fc{E}{X} Y  }
	}
	where
	\envMCm
	{
		\Fc{E}{X}
	}
	{
		\frc[s]{(1-X^2)}{X^2}
	}
	\envM
	{
		\Fc{O}{X}
	}
	{
		\frc[s]{(1-X^2)}{X^3}
	}
	and
	\envMCm
	{
		\Fc{S}{X}
	}
	{
		\frc[s]{(1-X^4)(1-X^6)}{X^{12}}	
	}
	and	
	$\Fc{S}{X}$
	is the generating series of the dimensions of the vector spaces of cusp forms on the full modular group.	
Specific examples for $r=2$ 		
	are given in \cite{GKZ06}	
	and
	a modern formulation is discussed in \cite{Brown21}
	(see also \cite{Tasaka16}).
However
	our computational results suggest the following  
	when we adopt $\setF[2]$ as the scalar field instead of $\setQ$:	%
	{(i)}
	the influence of modular forms disappears;
	but		
	{(ii)}		
	the Hoffman symbols $\bztO{\bfK}$ $(\bfK \in \midsHO{k})$ remain as basis elements 	
	with a Pascal triangle pattern.

It should  be noted that 
	the Broadhurst-Kreimer conjecture has 
	two equivalent formulations of vector and algebra (see \cite[Appendix]{IKZ06}).
The equivalence 
	requires 
	$\setQ[\ztO{2}]$ 
	is 
	isomorphic to the polynomial ring in one variable over $\setQ$.		
The isomorphy  
	does not hold 
	when $\setF[2]$ is the scalar field 
	as will be mentioned in the final section,
	and
	we will consider only the vector formulation in this paper. 	

It should also be noted that
	\etalTx{Bl\"umlein} \cite{BBV10}
	provided a data mine 
	for not only MZVs but also Euler sums
	by experiments to Broadhurst-Kreimer type conjectures,
	in which
	it was verified that
	the union of EDS and duality relations 
	suffices to reduce the number of generators of  $\spZO{k}$ to $\cdmZO{k}$
	up to weight $22$:
	it was also verified up to $24$ by using modular arithmetic,
	and
	up to $26$ and more 
	with an additional conjecture and limited depths.	
The duality relations,
	which 
	are 	
	obtained by
	the integral definition of MZVs 
	and
	a change of variables,
	are 
	very useful to compute 
	because they can bring down the size of relations by about half. 
It has not been proved yet that
	the EDS relations include the duality relations,
	although 
	the inclusion is expected to be true conjecturally: 	
	in other words,
	we have not succeeded in understanding the duality of MZVs algebraically. 
The experimental approaches of \cite{BBV10} and ours 
	differ
	in the use of the duality relations.
	
The organization of this paper is as follows.
In \refSect{sectTwo}, 
	we 
	state
	exact definitions of the binary MZVs $\bzt{\bfK}$,
	the formal multiple zeta spaces $\bspZO{k}$
	and
	the quotient spaces $\bspZqT{k}{r}$.
We report 
	our computational results in \refSect{sectThree},	
	and
	explain 
	how our computer programs produce the results in \refSect{sectFour}.		
The programs are available at the open-source site GitHub.\footnote{%
\siteMyCode
}
\refSect{sectFive} 	
	is devoted to
	problems about formal multiple zeta spaces 	
	which arise from the computational results.	
In Appendix,	
	we describe
	an essential algorithm in our experiments,
	which employs 
	a conflict based search 
	and
	speeds up the Gaussian forward elimination 
	under certain conditions.

The computer only assists us in showing \refPrp{4_PRP1} by Gaussian elimination. 		
The dimension conjecture \refEq{1_PL_CnjDimW} is true 
 	if	
	we can theoretically show \refPrp{4_PRP1} for all weights $k$.	

\section{Formal multiple zeta space over $\setF[2]$}\label{sectTwo}
The formal multiple zeta space $\bspZO{k}$ of weight $k$ 
	is briefly 
	defined by
	\envHLineCm[2_PL_DefBspZ]
	{
		\bspZO{k}
	}
	{
		\frc{\text{  \{binary EDS relations\} }  } 
		     {\SpT{\bZtO{\bfK}}{ \bfK \in \midsO{k} }_{\setF[2]} }
	}
	where 
	$\bZtO{\bfK}$ are indeterminates.
That is,
	$\bspZO{k}$ is
	an $\setF[2]$-vector space generated by	
	formal symbols $\bztO{\bfK} \equiv \bZtO{\bfK}$ 
	that satisfy 
	binary variations of the EDS relations. 	
Eight equivalent statements 	
	are given
	in \cite[Theorem 2]{IKZ06} for the EDS relations.	 
In this paper,
	we choice 
	the statement (v) in the theorem 
	because
	the relations are all $\setZ$-linear and fewer in number.	

To define \refEq{2_PL_DefBspZ} exactly,	
	we require   
	the algebraic setup by Hoffman \cite{Hoffman97}	
	which 
	allows us the steady handling of
	two products, 
	the shuffle $\shs$
	and
	stuffle $\shh$:
	the latter is also called harmonic or quasi-shuffle.
Let $\spHN$
	be the polynomial ring $\setQ\bkA{x,y}$ in the two non-commutative variables $x$ and $y$.
We call each variable a letter,
	and
	a monomial in the variables a word.
The shuffle product $\shs$ is 
	a $\setQ$-bilinear product on $\spHN$,
	which satisfies 
	\envMTh
	{
		w
	}
	{ 
		w \shS 1 
	}
	{
		1 \shS w 
	} 
	and
	\envHLine[2_PL_DefShS]
	{
		a u \shS b v 
	}
	{
		a (u \shS b v ) + b (a u \shS v )
	}
	for any words $u,v,w \in \spHN$
	and
	letters $a,b \in \SetO{x,y}$.
Let $z_k$ denote
	a word $x^{k-1} y$ for any $k \geq1$,
	and
	let 
	$\spHrN$ be 
	the polynomial ring $\setQ\bkA{z_1,z_2,\ldots}$,
	or equivalently,
	the subring $\setQ + \spHN y $ in $\spHN$.
The stuffle product $\shh$ is 
	a $\setQ$-bilinear product on $\spHrN$,
	which satisfies 
	\envMTh
	{
		w
	}
	{ 
		w \shH 1 
	}
	{
		1 \shH w 
	} 
	and
	\envHLine[2_PL_DefShH]
	{
		z_i u \shH z_j v 
	}
	{
		z_i (u \shH z_j v ) + z_j (z_i u \shH v ) + z_{i+j} (u \shH v)
	}
	for any words $u,v,w \in \spHN$
	and
	integers $i,j \geq 1$. 
By induction on the lengths of words,
	both products are commutative and associative,
	and
	both $\spHrN[\shs]=(\spHrN,\shs)$ and $\spHrN[\shh]=(\spHrN,\shh)$ 
	are 
	commutative $\setQ$-algebras.		
We notice 	
	$\spHN[\shs]=(\spHN,\shs)$ is a parent space of $\spHrN[\shs]$.
Let $\spHrrN=\setQ + x \spHN y=\SpT{z_{\bfK}}{ \bfK \in \midsN }_{\setQ}$,		
	where $z_{\bfK}=z_{k_1} \cdots z_{k_r}$ and $z_{\midE}=1$.
Both $\spHrrN[\shs]=(\spHrrN,\shs)$ and $\spHrrN[\shh]=(\spHrrN,\shh)$ 		
	are subalgebras	
	since 
	$\spHrrN$
	is closed under $\shs$ and $\shh$.
The pair $(\spHrN,\spHrrN)$ of spaces 
	satisfies  
	the polynomial ring property in one variable:		
	the former 	
	is freely generated by $y$ over the latter 
	on each of $\shs$ and $\shh$.
We thus have
	\envPLinePd[2_PL_IsoSpH]
	{
		\spHrN[\shs]	\lnP{\simeq}	\spHrrN[\shs][y]
	,\qquad
		\spHrN[\shh]	\lnP{\simeq}	\spHrrN[\shh][y]
	}
See \cite{Reutenauer93} and \cite{Hoffman97}
	for proofs of \refEq{2_PL_IsoSpH},
	respectively.

We introduce the EDS relations stated in \cite[Theorem 2(v)]{IKZ06}.
Let $\hmsRO[]{}$ denote  
	a homomorphism from $\spHrN[\shs]$ to $\spHrrN[\shs]$,
	which 
	is defined by
	taking the constant term with respect to $y$ 
	in the first isomorphism of \refEq{2_PL_IsoSpH}:\footnote{
The homomorphism $\hmhRO[]{}$ of stuffle type
	exists as well,
	but 
	it is intractable  
	because
	EDS relations of that type 		
	are not always $\setZ$-linear:
	see \cite{IKZ06} (or \cite{Bachmann20,Kaneko19}) for details.  
}
	\envPLinePd[2_PL_DefReg]
	{
		\hmsRO[]{} :
		\spHrN[\shs] \lnP{\ni} w \lnP{=} \SmT{i=0}{m} w_i \shS y^{\shs i} 
		\quad\mapsto\quad
		w_0 \lnP{\in} \spHrrN
	}
Let 	
	$\emidsO{k} = \midsO{k} \cup \SetO{(\usbTx{k}{1,\ldots,1})}$,
	and
	let
	\envHLinePd
	{
		\epidsO{k}
	}
	{
		\Un{i,j \geq 0 \atop (i+j=k)} \emidsO{i} \times \midsO{j}
	}
For any pair $(\bfK,\bfL)$
	of mult-indices in $\epidsO{k}$,
	we define 	
	\envHLineDefPd[2_PL_DefEDS]
	{
		\rlDO{\bfK,\bfL}
	}
	{
		\hmsRO{ z_{\bfK} \shH z_{\bfL} } - \hmsRO{ z_{\bfK} \shS z_{\bfL} }
		\lnP{\in}
		\spHrrN
	}
The objective EDS relations of weight $k$ 
	are
	stated as
	\envHLineCm[2_PL_EqEDS]
	{
		\evZO{ \rlDO{\bfK,\bfL} }
	}
	{
		0
		\qquad
		( (\bfK,\bfL) \in \epidsO{k} )
	}
	where
	$\evZO[]{}:\spHrrN\to\setR$
	is
	the $\setQ$-linear map (or \initWd{evaluation map})  
	defined by 
	$\evZO{z_{\bfK}} = \zt{\bfK }$ $(\bfK \in \midsN)$.	
We have by \refEq{2_PL_DefReg}
	\envHLineCFPd
	{
		\hmsRO{w}
	}
	{
		w
		\qquad
		(w\in\spHrrN)
	}
	{
		\hmsRO{ y^m \shS z_{\bfM} }
	}
	{
		0
		\qquad
		(m>0, \bfM \in \midsN)
	}	
 We can thus divide \refEq{2_PL_EqEDS} into two parts:
	\envHLineCFCmNme
	{\label{2_PL_EqEDSfin}
		\evZO{ z_{\bfK} \shH z_{\bfL} } -  \evZO{ z_{\bfK} \shS z_{\bfL} }
	}
	{
		0
		\qquad
		( (\bfK,\bfL) \in \pidsO{k} )	
	}
	{\label{2_PL_EqEDSprp}
		\evZO{ \hmsRO{ y^m \shH z_{\bfM} } }
	}
	{
		0
		\qquad
		( 0 < m < k-1, \bfM \in \midsO{k-m} )
	}
	where
	\envMPd{
		\pidsO{k}
	}{
		\Un[t]{i,j \geq 0 \atop (i+j=k)} \midsO{i} \times \midsO{j}	
	}	
The relations in \refEq{2_PL_EqEDSfin} 
	are
	called 
	the finite double shuffle (FDS) relations,
	because
	MZVs are defined by  
	\envM
	{
		\zt{k_1, \ldots, k_r}
	}
	{
		\Sm{m_1>\cdots>m_r>0} \frc[s]{ m_1^{k_1} \cdots m_r^{k_r} }{1}		
	} 
	and	
	finite (or convergent) at $\bfK \in \midsN$.	
The FDS relations do not suffice to give all relations of MZVs.
For instance,
	we can not obtain any relation in weight $3$,
	in particular,
	the simplest formula $\zt{2,1}=\zt{3}$. 	
Therefore
	the relations in \refEq{2_PL_EqEDSprp}	
	are 
	essential to the EDS conjecture.		

A little more notions 
	are required for \refEq{2_PL_DefBspZ},
	which 
	are analogs of the notions mentioned above in $\setZ$-module and $\setF[2]$-vector.
Let $\ispHN$
	denote the subring $\setZ\bkA{x,y}$ in $\setQ\bkA{x,y}$.
We set  
	\envPLineCm
	{	
		\ispHrrN	\lnP{=}		\SpT{z_{\bfK}}{ \bfK \in \midsN }_{\setZ}	
	,\qquad
		\bspUN	\lnP{=}	\SpT{\bZtO{\bfK}}{ \bfK \in \midsN }_{\setF[2]}
	}
	to define 
	a canonical map from $\ispHrrN$ to $\bspUN$
	which is induced by modulo $2$ arithmetic:
	\envPLinePd[2_PL_DefLin2]
	{
		\lnCO[]{} :
		\ispHrrN \lnP{\ni} w \lnP{=} \Sm{\bfK \in \midsN} c_{\bfK} z_{\bfK} 
		\quad\mapsto\quad
		\Sm{\bfK \in \midsN} (c_{\bfK}\,\modd 2) \bZt{\bfK}  \lnP{\in} \bspUN
	}
For any pair $(\bfK,\bfL) \in \pidsO{k}$,
	the elements $z_{\bfK} \shH z_{\bfL}$ and $z_{\bfK} \shS z_{\bfL}$ belong to 
	\envMOCm
	{
		\ispHrrN
	}
	and
	the element $\lnCO{\rlDO{\bfK,\bfL}}$ 
	is well-defined.
For $0 < m < k-1$ and $\bfM \in \midsO{k-m}$,
	$y^m \shH z_{\bfM}$ belongs  to 
	\envMOCm
	{
		\SpT{y^n z_{\bfN} }{ n \geq 0, \bfN \in \midsN }_\setZ
	}
	and
	$\lnCO{ \hmsRO{ y^m \shH z_{\bfM} } }$ 
	is well-defined 
	if 
	\envHLineCmPt{\in}	
	{
		\hmsRO{ y^n z_{\bfN} }
	}
	{
		\ispHrrN	\qquad	( n > 0, \bfN \in \midsN )
	}
	which 
	holds by \cite[Proposition 8]{IKZ06} (see \refEq{4_PL_EqHmsR} below).		
Consequently, 		
	\envHLineDef
	{
		\bspEO{k} 
	}
	{ 
		\SpT{ \lnCO{\rlDO{\bfK,\bfL}} }{ (\bfK,\bfL) \in \epidsO{k} }_{\setF[2]}
	\lnP{\subset}
		\bspUO{k}
	}	
	is well-defined.
	
We are in a position to define 
	\refEq{2_PL_DefBspZ}.	
	
\begin{definition}\label{2_DEF1}
For a weight $k$,
	we define 
	the formal multiple zeta space 	
	by 
	\envHLineDefPd[2_DEF1_FMZsp]
	{
		\bspZO{k}
	}
	{
		\bspUO{k} / \bspEO{k}
	}
For a mult-index $\bfK\in \midsO{k}$,
	we denote by  
	$\bzt{\bfK}$ 
	the element in $\bspZO{k}$ 
	which is congruent to $\bZt{\bfK}$ modulo $\bspEO{k}$.
We call 
	$\bzt{\bfK}$ a binary multiple zeta symbol 
	or 
	a binary MZS.
\end{definition}

Let 
	$\bevHO[]{}$	
	denote
	the natural homomorphism from $\OPs{k\geq0} \bspUO{k}$ to $\OPs{k\geq0} \bspZO{k}$:
	each component 	
	is 	
	the canonical map of \refEq{2_DEF1_FMZsp}.
We define
	the \initWd{binary evaluation map}
	by 
	\envMPd{
		\bevZO[]{} 
	}{
		\bevHO[]{} \circ \lnCO[]{}
	}	
The binary EDS relations of weight $k$ 
	are then stated as 	
	\envHLinePd[2_PL_EqBEDS]
	{
		\bevZO{ \rlDO{\bfK,\bfL} }
	}
	{
		0
		\qquad
		( (\bfK,\bfL) \in \epidsO{k} )
	}
We list 	
	some examples of the original and binary EDS relations for weights $k\leq 4$  	
	in \refTab{2_Tbl_EDS}.

Let $\bspZT{k}{r}$
	denote
	the vector subspace $\SpT{\bzt{\bfK}}{\bfK \in \midsO{k}, \dpt{\bfK} \leq r}_{\setF[2]}$ 	
	as introduced in the first section.	
We end this section 
	with 
	the definition of the graded pieces satisfying 
	the direct sum decomposition \refEq{1_PL_CgDSDZ}.
\begin{definition}\label{2_DEF2}
For a weight $k$,	 
	we define 
	the depth graded formal multiple zeta spaces 
	by
	\envHLineCmDef[2_DEF2_qFMZsp]
	{
		\bspZqT{k}{r}
	}
	{
		\bspZT{k}{r} / \bspZT{k}{r-1}
	\qquad
		(k >  r \geq 0)
	}	
	where
	$\bspZT{k}{-1} = \SetO{0}$.
\end{definition}

\begin{table}[!t]\renewcommand{\arraystretch}{1.4}	
\begin{center}{
	\caption{EDS relations in $\spZO{k}$ and $\bspZO{k}$ for weights $k\leq4$.}\label{2_Tbl_EDS}
	\mbox{}\\\vPack[5]\begin{tabular}{|c|c|c|}  \hline
		$\bfK,\bfL$ &Original EDS relation (over $\setZ$)&Binary EDS relation	 (over $\setF[2]$)		\\\hline  	
		$(1),(2)$     &$-\zt{2,1}+\zt{3}=0$	& $\bzt{2,1}+\bzt{3}=0$	\\\hline
		$(1),(3)$     &$-\zt{2,2}-\zt{3,1}+\zt{4}=0$	& $\bzt{2,2}+\bzt{3,1}+\bzt{4}=0$	\\\hline
		$(1),(2,1)$     &$-\zt{2,1,1}+\zt{2,2}+\zt{3,1}=0$	& $\bzt{2,1,1}+\bzt{2,2}+\bzt{3,1}=0$	\\\hline
		$(1,1),(2)$     &$\zt{2,1,1}-\zt{2,2}-\zt{3,1}=0$	& $\bzt{2,1,1}+\bzt{2,2}+\bzt{3,1}=0$	\\\hline
		$(2),(2)$     &$-4\zt{3,1}+\zt{4}=0$	& $\bzt{4}=0$	\\\hline
	\end{tabular} 
}\end{center}
\end{table}

\section{Computational result} \label{sectThree}
We report our computational results.
How we obtain them
	will be explained in the next section.

We begin with  
	a typical result related to \refEq{1_PL_CnjDimW}.		
\begin{experiment}\label{3_EXP1}
For any weight $k$ with $2 \leq k \leq 22$, 	
	we verify 
	$\bztO{\midsHO{k}}$ is a basis of $\bspZO{k}$,
	and
	\envOTLinePd[3_EXP1_Eq]
	{
		\dimT{\setF[2]}{\bspZO{k}}
	}
	{
		2^{k-2} - \dimT{\setF[2]}{\bspEO{k}}
	}
	{
		\cdmZO{k}
	}	
\end{experiment}

The EDS conjecture 
	states that,
	for every weight $k$,
	the relations in \refEq{2_PL_EqEDS}
	suffice to reduce the number of generators of  $\spZO{k}$ to $\cdmZO{k}$:
	\envHLineCmPt[3_PL_CnjEDS]{\geq}
	{
		\dimT{\setQ}{\spEO{k}} 
	}
	{
		2^{k-2} - \cdmZO{k}
	}		
	where 
	\envMPd
	{
		\spEO{k} 
	}
	{ 
		\SpT{ \rlDO{\bfK,\bfL} }{ (\bfK,\bfL) \in \epidsO{k} }_{\setQ}
	}
This can be confirmed by 
	\refExp{3_EXP1},	
	as follows.	
We denote 
	by $\ispEO{k} = \SpT{ \rlDO{\bfK,\bfL} }{ (\bfK,\bfL) \in \epidsO{k} }_{\setZ}$
	the $\setZ$-module counterpart of $\spEO{k}$.
Since	
	$\setQ$ is the field of fractions of $\setZ$
	and
	$\lnCO[]{}$ is a surjective homomorphism from $\ispEO{k}$ to $\bspEO{k}$,	
	\envHLineCm
	{
		\dimT{\setQ}{\spEO{k}}
	}
	{
		\rnkT{\setZ}{\ispEO{k}}
	\lnP{\geq}
		\dimT{\setF[2]}{\bspEO{k}}
	}
	which,
	together with \refEq{3_EXP1_Eq},
	proves 
	\refEq{3_PL_CnjEDS} for $k \leq 22$.

We recall	
	$\cdmBZqT{k}{r}=\bnm{k-2r}{r}$
	that is
	the number of the Hoffman mult-indices of weight $k$ and depth $r$.
We define $\bspUT{k}{r} = \SpT{\bZtO{\bfK}}{\bfK \in \midsO{k}, \dpt{\bfK} \leq r}_{\setF[2]} \subset \bspUO{k}$,
	and	
	\envHLinePd
	{
		\bspEqT{k}{r}
	}
	{
		\bkR{ \bspUT{k}{r} \cap \bspEO{k} } / \bspUT{k}{r-1}
	}

	
The main result 
	is 
	a refinement of \refExp{3_EXP1}.	
Taking the sum for $r=1,\ldots,k-1$ in \refEq{3_EXP2_Eq}
	induces \refEq{3_EXP1_Eq}
	because of \refEq{1_PL_CgDSDZ}:
	note
	that $\bspZqT{k}{0} = \SetO{0}$ unless $k = 0$.

\begin{experiment}\label{3_EXP2}
For any weight $k$ and depth $r$ with $1 \leq r < k \leq 22$,
	we verify 
	$\bztqO{\midsHT{k}{r}}$ is a basis of $\bspZqT{k}{r}$,
	and
	\envOTLinePd[3_EXP2_Eq]
	{
		\dimT{\setF[2]}{\bspZqT{k}{r}}
	}
	{
		\bnm{r-1}{k-2} - \dimT{\setF[2]}{\bspEqT{k}{r}}
	}
	{
		\cdmBZqT{k}{r}
	}
\end{experiment}

The first equality in \refEq{3_EXP2_Eq}
	is
	by the isomorphism theorems.		
In fact,
	we have
	\envHLineFCmPt[3_PL_SEqBZSa]{\simeq}
	{	
		\bspZT{k}{r} / \bspZT{k}{r-1}
	}
	{
		\raMTx{4}{ \bkR[]{ \bspUT{k}{r} / \bkR{ \bspUT{k}{r} \cap \bspEO{k} } } }
		\mcsgsb[B]{/} 
		\raMTx{-4}{ \bkR[]{ \bspUT{k}{r-1} / \bkR{ \bspUT{k}{r-1} \cap \bspEO{k} } } }
	}
	{\rule{0pt}{15pt}
		\raMTx{4}{ \bspUT{k}{r}  }
		\mcsgsb[B]{/} 
		\raMTx{-4}{  \bkR{ \bspUT{k}{r-1} +  \bkR[]{ \bspUT{k}{r} \cap \bspEO{k} } } }
	}
	{\rule{0pt}{15pt}
		\raMTx{4}{ \bkR[]{ \bspUT{k}{r} / \bspUT{k}{r-1}  } }
		\mcsgsb[B]{/} 
		\raMTx{-4}{ \bkR[]{ \bkR{ \bspUT{k}{r}  \cap \bspEO{k} } / \bspUT{k}{r-1} } }
	}
	and
	\envHLinePdPt{\simeq}	
	{
		\bspZqT{k}{r}
	}
	{
		\raMTx{4}{ \bkR[]{ \bspUT{k}{r} / \bspUT{k}{r-1}  } }
		\mcsgsb[B]{/} 
		\raMTx{-4}{ \bspEqT{k}{r} }
	}
Since
	$\bnm{r-1}{k-2} = \nmSet{\bspUT{k}{r} / \bspUT{k}{r-1}}$ 
	by 
	counting the number of the mult-indices of weight $k$ and depth $r$,
	we obtain
	the desired equality.

We demonstrate the numbers $\cdmBZqT{k}{r}$ for $k\leq22$ in \refTab{3_Tbl_EXP12_dd}.
They 
	are expressed in terms of binomial coefficients,
	and
	we can observe a (shifted) Pascal triangle pattern:	
	the column $r=0$ has the sequence $(1)$ from the row $k=0$,
	the column $r=1$ has $(1,1)$ from $k=2$,
	the column $r=2$ has $(1,2,1)$ from $k=4$,
	the column $r=3$ has $(1,3,3,1)$ from $k=6$,
	and so on.
For comparison,	
	the dimensions of ${\spZqT{k}{r}}$ conjectured in \refEq{1_PL_EqCDmZqBK}	
	are listed in \refTab{3_Tbl_EXP12_dd_BK}.

\begin{table}[!t]\renewcommand{\arraystretch}{1.1} \arraycolsep=8.0pt
\begin{center}
	\caption{
		The numbers $\cdmBZqT{k}{r}$ for $0 \leq r < k \leq 22$:	%
		the unlisted numbers $\cdmBZqT{k}{r}$ $(r>11)$ are $0$.	
		The total number of each row is $\cdmZO{k}$
		and
		that of each column is $2^{r}$ (for $r\leq 7$).
	}\label{3_Tbl_EXP12_dd}
	\vWiden
	${\left.\begin{array}{|c|cccccccccccc|c|}\hline	
		k/r & 0 & 1 & 2 & 3 & 4 & 5 & 6 & 7 & 8 & 9 & 10 & 11 & \text{Total} \\\hline
		0 & 1 & 0 & 0 & 0 & 0 & 0 & 0 & 0 & 0 & 0 & 0 & 0 & 1 \\ 
		1 & 0 & 0 & 0 & 0 & 0 & 0 & 0 & 0 & 0 & 0 & 0 & 0 & 0 \\ 
		2 & 0 & 1 & 0 & 0 & 0 & 0 & 0 & 0 & 0 & 0 & 0 & 0 & 1 \\ 
		3 & 0 & 1 & 0 & 0 & 0 & 0 & 0 & 0 & 0 & 0 & 0 & 0 & 1 \\ 
		4 & 0 & 0 & 1 & 0 & 0 & 0 & 0 & 0 & 0 & 0 & 0 & 0 & 1  \\ 
		5 & 0 & 0 & 2 & 0 & 0 & 0 & 0 & 0 & 0 & 0 & 0 & 0 & 2 \\ 
		6 & 0 & 0 & 1 & 1 & 0 & 0 & 0 & 0 & 0 & 0 & 0 & 0 & 2 \\ 
		7 & 0 & 0 & 0 & 3 & 0 & 0 & 0 & 0 & 0 & 0 & 0 & 0 & 3 \\ 
		8 & 0 & 0 & 0 & 3 & 1 & 0 & 0 & 0 & 0 & 0 & 0 & 0 & 4 \\ 
		9 & 0 & 0 & 0 & 1 & 4 & 0 & 0 & 0 & 0 & 0 & 0 & 0 & 5 \\ 
		10 & 0 & 0 & 0 & 0 & 6 & 1 & 0 & 0 & 0 & 0 & 0 & 0 & 7 \\ 
		11 & 0 & 0 & 0 & 0 & 4 & 5 & 0 & 0 & 0 & 0 & 0 & 0 & 9 \\ 
		12 & 0 & 0 & 0 & 0 & 1 & 10 & 1 & 0 & 0 & 0 & 0 & 0 & 12 \\ 
		13 & 0 & 0 & 0 & 0 & 0 & 10 & 6 & 0 & 0 & 0 & 0 & 0 & 16 \\ 
		14 & 0 & 0 & 0 & 0 & 0 & 5 & 15 & 1 & 0 & 0 & 0 & 0 & 21 \\ 
		15 & 0 & 0 & 0 & 0 & 0 & 1 & 20 & 7 & 0 & 0 & 0 & 0 & 28 \\ 
		16 & 0 & 0 & 0 & 0 & 0 & 0 & 15 & 21 & 1 & 0 & 0 & 0 & 37 \\ 
		17 & 0 & 0 & 0 & 0 & 0 & 0 & 6 & 35 & 8 & 0 & 0 & 0 & 49 \\ 
		18 & 0 & 0 & 0 & 0 & 0 & 0 & 1 & 35 & 28 & 1 & 0 & 0 & 65 \\
		19 & 0 & 0 & 0 & 0 & 0 & 0 & 0 & 21 & 56 & 9 & 0 & 0 & 86 \\ 
		20 & 0 & 0 & 0 & 0 & 0 & 0 & 0 & 7 & 70 & 36 & 1 & 0 & 114 \\ 
		21 & 0 & 0 & 0 & 0 & 0 & 0 & 0 & 1 & 56 & 84 & 10 & 0 & 151 \\ 
		22 & 0 & 0 & 0 & 0 & 0 & 0 & 0 & 0& 28 & 126 & 45 & 1 & 200 \\ \hline
		\text{Total} & 1 & 2 & 4 & 8 & 16 & 32 & 64 & 128 & - & - & - & - & \multicolumn{1}{c}{} \\ \cline{1-13}
	\end{array}\right.}$
\end{center}
\end{table}
\begin{table}[!t]\renewcommand{\arraystretch}{1.1} \arraycolsep=8.0pt
\begin{center}
	\caption{
		The conjectural numbers $\dimO{\spZqT{k}{r}}$ for $0 \leq r < k \leq 22$:	
		the unlisted numbers $\dimO{\spZqT{k}{r}}$ $(r > 11)$ are $0$.	
	}\label{3_Tbl_EXP12_dd_BK}
	\vWiden
	${\left.\begin{array}{|c|cccccccccccc|c|}\hline	
		k/r & 0 & 1 & 2 & 3 & 4 & 5 & 6 & 7 & 8 & 9 & 10 & 11 & \text{Total} \\\hline
		0 & 1 & 0 & 0 & 0 & 0 & 0 & 0 & 0 & 0 & 0 & 0 & 0 & 1 \\ 
		1 & 0 & 0 & 0 & 0 & 0 & 0 & 0 & 0 & 0 & 0 & 0 & 0 & 0 \\ 
		2 & 0 & 1 & 0 & 0 & 0 & 0 & 0 & 0 & 0 & 0 & 0 & 0 & 1 \\ 
		3 & 0 & 1 & 0 & 0 & 0 & 0 & 0 & 0 & 0 & 0 & 0 & 0 & 1 \\ 
		4 & 0 & 1 & 0 & 0 & 0 & 0 & 0 & 0 & 0 & 0 & 0 & 0 & 1  \\ 
		5 & 0 & 1 & 1 & 0 & 0 & 0 & 0 & 0 & 0 & 0 & 0 & 0 & 2 \\ 
		6 & 0 & 1 & 1 & 0 & 0 & 0 & 0 & 0 & 0 & 0 & 0 & 0 & 2 \\ 
		7 & 0 & 1 & 2 & 0 & 0 & 0 & 0 & 0 & 0 & 0 & 0 & 0 & 3 \\ 
		8 & 0 & 1 & 2 & 1 & 0 & 0 & 0 & 0 & 0 & 0 & 0 & 0 & 4 \\ 
		9 & 0 & 1 & 3 & 1 & 0 & 0 & 0 & 0 & 0 & 0 & 0 & 0 & 5 \\ 
		10 & 0 & 1 & 3 & 3 & 0 & 0 & 0 & 0 & 0 & 0 & 0 & 0 & 7 \\ 
		11 & 0 & 1 & 4 & 3 & 1 & 0 & 0 & 0 & 0 & 0 & 0 & 0 & 9 \\ 
		12 & 0 & 1 & 3 & 6 & 2 & 0 & 0 & 0 & 0 & 0 & 0 & 0 & 12 \\ 
		13 & 0 & 1 & 5 & 6 & 4 & 0 & 0 & 0 & 0 & 0 & 0 & 0 & 16 \\ 
		14 & 0 & 1 & 5 & 9 & 4 & 2 & 0 & 0 & 0 & 0 & 0 & 0 & 21 \\ 
		15 & 0 & 1 & 6 & 8 & 10 & 3 & 0 & 0 & 0 & 0 & 0 & 0 & 28 \\ 
		16 & 0 & 1 & 5 & 14 & 11 & 6 & 0 & 0 & 0 & 0 & 0 & 0 & 37 \\ 
		17 & 0 & 1 & 7 & 13 & 18 & 7 & 3 & 0 & 0 & 0 & 0 & 0 & 49 \\ 
		18 & 0 & 1 & 6 & 19 & 18 & 17 & 4 & 0 & 0 & 0 & 0 & 0 & 65 \\
		19 & 0 & 1 & 8 & 17 & 31 & 19 & 10 & 0 & 0 & 0 & 0 & 0 & 86 \\ 
		20 & 0 & 1 & 7 & 25 & 30 & 35 & 12 & 4 & 0 & 0 & 0 & 0 & 114 \\ 
		21 & 0 & 1 & 9 & 22 & 48 & 37 & 29 & 5 & 0 & 0 & 0 & 0 & 151 \\ 
		22 & 0 & 1 & 8 & 32 & 45 & 65 & 33 & 16& 0 & 0 & 0 & 0 & 200 \\ \hline
	\end{array}\right.}$
\end{center}
\end{table}

Refinements	
	of the EDS conjecture 		
	have been proposed.	
\etalTx{Minh} \cite{MJPO00}
	conjectured  	
	that
	a part of the EDS relations obtained from 
	\envHLine[3_PL_DefPidsMJPO]
	{
		\epidsMO{k}
	}
	{
		\pidsO{k} \cup (\emidsO{1} \times \midsO{k-1})
	}
	is 
	a right candidate, 	
	and 
	verified it up to $k=16$.
The relations 
	\envHLine
	{
		\evZO{ \rlDO{\bfK,\bfL} }
	}
	{
		0
		\qquad
		( (\bfK,\bfL) \in \emidsO{1} \times \midsO{k-1})
	}
	are 
	known as Hoffman's relations (\cite{Hoffman97}),	
	and
	their conjecture says that 
	FDS relations and Hoffman's relations suffice to  give all relations among MZVs.
\etalTx{Kaneko} \cite{KNT08} 
	conjectured  
	the above relations are too much,
	i.e.,
	a smaller part obtained from 
	\envHLine[3_PL_DefPidsKNT]
	{
		\epidsKO{k}
	}
	{
		(\SetO{(3),(2,1)} \times \midsO{k-3})
		\cup
		(\SetO{(2)} \times \midsO{k-2})
		\cup
		(\emidsO{1} \times \midsO{k-1})
	}
	is 
	a right candidate.
They
	verified it up to $k=20$.

In the space $\bspZO{k}$,		
	neither 
	the relations obtained from \refEq{3_PL_DefPidsMJPO} 
	nor 
	those obtained from \refEq{3_PL_DefPidsKNT}
	suffice 
	to give all relations among binary MZSs.
	
\begin{experiment}\label{3_EXP3}
Let $\bullet \in \SetO{\letKNT,\letMJPO}$
	and
	let 
	$\bspEBO{k}=\SpT{ \lnCO{\rlDO{\bfK,\bfL}} }{ (\bfK,\bfL) \in \epidsBO{k} }_{\setF[2]}$.
There exist weights $k \leq 22$ 
	such that
	\envHLinePdPt[3_EXP3_Eq]{<}
	{
		\dimT{\setF[2]}{\bspEBO{k}}
	}
	{
		2^{k-2} - \cdmZO{k}
	}
\end{experiment}
\begin{table}[!t]\renewcommand{\arraystretch}{1.1} \arraycolsep=20pt
\begin{center}
	\caption{
		The numbers $\cdmBZBO{k}$ $(\bullet \in \SetO{\letKNT,\letMJPO})$
		with
		$\cdmZO{k}$:
		they are same when $k \leq 6$.
	}\label{3_Tbl_EXP3}
	\vWiden
	$\left.\begin{array}{|c|ccc|}\hline 
		k \rule[-6pt]{0pt}{20pt}& \cdmBZKO{k} & \cdmBZMO{k} & \cdmZO{k} \\\hline
		%
		7 & 4 & 4  & 3 \\
		8 & 6 & 4  & 4 \\
		9 & 8 & 6  & 5 \\
		10 & 12 & 8  & 7 \\
		11 & 21 & 10  & 9 \\
		12 & 30 & 14  & 12 \\
		13 & 44 & 18  & 16 \\
		14 & 66 & 24  & 21 \\
		15 & 100 & 33  & 28 \\
		16 & 140 & 42  & 37 \\
		17 & 208 & 57  & 49 \\
		18 & 300 & 75  & 65 \\ 
		19 & 441  & 99  & 86 \\ 
		20 & 644  & 132  & 114 \\ 
		21 & -  & 174  & 151 \\ 
		22 & -  & 231  & 200 \\\hline 
	\end{array}\right.$
\end{center}
\end{table}

Computational results of $\cdmBZBO{k} = 2^{k-2}-\dimO{\bspEBO{k}}$ 
	are 
	shown in \refTab{3_Tbl_EXP3}.
In general,
	\envOTLineThPdPt{>}{
		\cdmBZKO{k}
	}{ 
		\cdmBZMO{k} 
	}{ 
		\cdmZO{k}
	}
We can find that	
	the sequence $\bkR{\cdmBZMO{k}}_{0 \leq k \leq 22}$ 
	has 
	a quasi Fibonacci-like rule,
	\envHLineCm[3_PL_EqFibMJPO]
	{
		\cdmBZMO{k} 
	}
	{ 
		\cdmBZMO{k-2} + \cdmBZMO{k-3} + \del_{M,k}	
	}
	where	
	$M=\SetO{7,15}$
	and
	$\del_{M,k}$ is the Kronecker delta function 
	defined by $\del_{M,k} = 1$ if $k \in M$ and $\del_{M,k} = 0$ otherwise.		
It appears 
	that 
	$\bkR{\cdmBZKO{k}}_{0 \leq k \leq 22}$ 
	does not have an obvious law.

\section{Computer program} \label{sectFour}
Our computer programs,
	that perform   
	the Gaussian forward elimination on the linear combinations in $\bspEO{k}$,
	show the following proposition.
\begin{proposition}\label{4_PRP1}
Let $k$ and $r$ be a weight and depth,
	respectively,	
	with 
	$r < k \leq 22$.
For a mult-index $\bfK$ in $\midsT{k}{r}$,		
	the following statements hold.
\mbox{}\vWiden[3]\mbox{}\\{\bf (i)}
If $\bfK \notin \midsHT{k}{r}$,	
	there exists a combination $c \in \bspUT{k}{r} \cap\bspEO{k}$
	such that
	\envHLinePdPt[4_PRP1_IncBHS]{\in}
	{
		\bZtO{\bfK} 
	}
	{
		c 
		+ 	
		\SpT{ \bZtO{\bfH} }{ \bfH \in \midsHT{k}{r} \cup \midsHT{k}{r-1} \cup \cdots \cup \midsHT{k}{\bkF{\frc[s]{3}{k}}} }_{\setF[2]} 
	}	
\mbox{}{\bf (ii)}
If $\bfK \in \midsHT{k}{r}$,
	there exists no combination $c$ such as $\mathrm{\refEq{4_PRP1_IncBHS}}$.
\mbox{}\vWiden[3]\par
Here
	$\bkF{\cdot}$ is the floor function	
	defined by $\bkF{t} = \MaxT{ a \in \setZ }{a \leq t}$ for a real number $t$.	
\end{proposition}

\refPrp{4_PRP1} 
	verifies 	
	\refExp{3_EXP2}.
Suppose $\bfK \in \midsT{k}{r} \setminus \midsHT{k}{r}$.
By the statement {\bf (i)},		
	\envHLineCmPt[4_PL_IncBZS]{\in}
	{
		\bztO{\bfK} 
	}
	{
		\SpT{ \bztO{\bfH} }{ \bfH \in \midsHT{k}{r} \cup \midsHT{k}{r-1} \cup \cdots \cup \midsHT{k}{\bkF{\frc[s]{3}{k}}} }_{\setF[2]} 
	}
	or
	\envHLineCmPt[4_PL_IncBZSq]{\in}
	{
		\bztqO{\bfK} 
	}
	{
		\SpT{ \bztqO{\bfH} }{ \bfH \in \midsHT{k}{r} }_{\setF[2]} 
	}		
	which,
	together with the statement {\bf (ii)},		
	implies
	$\bztqO{\midsHT{k}{r} }$ 
	is 
	a basis of $\bspZqT{k}{r}$
	for $r < k \leq 22$.

Imaginarily, 		
	the Gaussian elimination 
	can 
	elucidate
	any vector space whose 
	corresponding matrix (or set of defining linear combinations) is clearly given:
	but
	practically,
	it is limited to 
	a space that are not too big. 
The bound $k=22$ in \refPrp{4_PRP1} 
	indicates 
	a performance threshold of our computing environments.
Below
	we will describe the environments and prove \refPrp{4_PRP1}.

The programs	
	are 
	written almost by 
	Python language
	and partly by
	Cython language.
The machine is as follows:
	a Linux-based PC having
	two CPUs with $12$-core at 2.70GHz (Intel Xeon Gold 6226)
	and
	a $3$TB RAM.
The package of the programs
	is available at \siteMyCode.	

\begin{figure}[!t]
\begin{center}
	\caption{
		Layout of our package for the executable files.
	}\label{3_Fig_LayPac}
	\vWiden
	\includegraphics[keepaspectratio, scale=0.5]{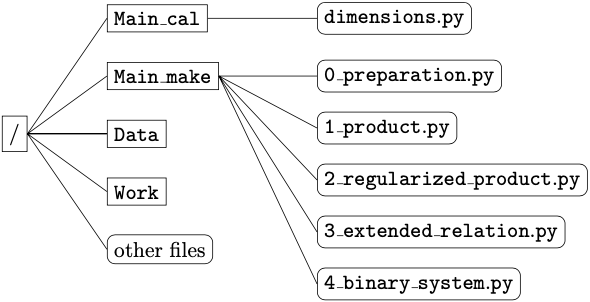}
\end{center}
\end{figure}

The executable files
	are in the directories 
	named as
	$\dirM$ 
	and 
	$\dirC$.
The former contains five files 
	that produce 
	datas of binary systems (or binary matrices) 
	obtained from the binary EDS relations,
	and
	the latter contains one file 
	that calculates 
	dimensions of $\bspZO{k}$ and $\bspZqT{k}{r}$
	(or row echelon forms of the corresponding binary matrices).
The produced datas	
	are stocked in $\dirD$,
	almost of which are saved in Python pickle format to reduce data size. 
Class files
	in which essential precesses are performed
	are stored in $\dirW$.
Files of config, license and readme are also placed
	in the root directory of the package. 
(See \refFig{3_Fig_LayPac} for a layout of the package).

We have a convenient expression for a linear combination in $\bspZO{k}$
	since		
	$\setF[2]$ consists of only two elements.	
A subset $\bfJ[C]$ in $\midsO{k}$
	is identified with 
	a combination such as
	\envPLinePd[4_PL_One2One]
	{
		\bfJ[C] 	
	\quad\lnP{\longleftrightarrow}\quad
		\Sm{\bfK \in \bfJ[C]} \bzt{\bfK} 		
	}
For instance,
	$\bfJ[C]_1 = \SetO{(2,1,1),(2,2),(3,1)}$
	corresponds to
	$\bzt{2,1,1}+\bzt{2,2}+\bzt{3,1}$
	and
	$\bfJ[C]_2 = \SetO{(2,2),(3,1),(4)}$
	corresponds to
	$\bzt{2,2}+\bzt{3,1}+\bzt{4}$.
By \refEq{4_PL_One2One},	
	the symmetric difference $\soSD$ of two sets
	is 
	equivalent to the plus of two combinations:	
	\envM
	{
		\bfJ[C]_1 \soSD \bfJ[C]_2 
	}
	{ 
		(\bfJ[C]_1 \setminus \bfJ[C]_2) \cup (\bfJ[C]_2 \setminus \bfJ[C]_1) = \SetO{(2,1,1),(4)}
	}
	corresponds to 		
	\envMPd
	{
		\bkR{\bzt{2,2}+\bzt{3,1}+\bzt{4}}+\bkR{\bzt{2,1,1}+\bzt{2,2}+\bzt{3,1}}
	}	
	{
		\bzt{2,1,1}+\bzt{4}
	}
The expression \refEq{4_PL_One2One}
	is also applied to
	a linear relation
	$\Sm{\bfK \in \bfJ[C]} \bzt{\bfK} 	= 0$ in the same way.
We compute 
	binary EDS relations (or defining combinations in $\bspEO{k}$)
	through \refEq{4_PL_One2One}
	with 
	the set datatype in Python.
This set based expression 	
	can be
	realized by built-in objects.\footnote{
	We use
	\emph{frozenset} 
	and 
	\emph{s.symmetric\_difference(t)} (or the operator notation `\emph{s}\textasciicircum\emph{t}'),		
	where
	\emph{frozenset} is an immutable datatype for set datas and \emph{s,t} are its instances.
}

We will explain  
	the executable files
	and
	report their statics.	
We do not mention 
	actual command lines to use the files in a linux OS,
	but
	we can find them 
	in the beginning of each file. 

\subsection{Executable file in $\dirM$}
We will require many maps to save midway datas 	
	for the binary linear systems induced from the binary EDS relations.		
The prime reason is 
	that,
	by \refEq{2_PL_DefEDS},
	each EDS relation
	is 	
	composed of a combination of $\shS[]$, $\shH[]$ and $\hmsRO[]{}$. 	
For the maps or the midway datas,
	we will use dictionary datatype,
	which 	
	is a built-in object in Python 
	and
	consists of a collection of 
	tuples of two objects called `key' and `value':
	a key-object is mapped to its associated value-object.		

The file $\fleM[0]$ prepares 
	two dictionary datas for each weight $k \leq 22$.
Let $\bkS{n}$ 
	denote the set $\SetO{1,\ldots,n}$ 	
	for a positive integer $n$. 		
One data
 	gives 
	a one-to-one mapping 
	from
	the integers in $\bkS{2^k}$ to the words of degree $k$,		
	and
	another data gives 
	a one-to-one mapping from
	the integers in $\bkS{2^{k-1}}$ to the mult-indices in $\UnT[t]{r=1}{k} \setN^r$ of weight $k$:
	if $n\in\bkS{2^{k-1}}$ 
	and 
	the associated mult-index is $\bfK$,	
	the associated word is $z_{\bfK}$.
The objects of the set 
	which  	
	our programs select for the set based expression in the left of \refEq{4_PL_One2One} 	
	are	
	the integers (for which the integer datatype is necessary) 
	instead of the mult-indices and words (for which the tuple and string datatypes are necessary),
	because
	the integer datatype is reasonable in data size and running time.		

The file $\fleM[1]$ 
	creates 
	dictionary datas for shuffle and stuffle products.
The defining equations \refEq{2_PL_DefShS} and \refEq{2_PL_DefShH}
	suggest
	that
	creating datas of shuffle  
	will take 
	more time 
	since 
	shuffle products can contain more terms.
For a speed-up, 	
	we improve \refEq{2_PL_DefShS}:
	\envLineCm[4_PL_EqShS]
	{
		a_1 \cdots a_m \shS b_1 \cdots b_n 
	}
	{
		\tpSm{i+j=l}{0 \leq i \leq \Min{l,m} \atop 0 \leq j \leq \Min{l,n} }
		\bkR[a]{  a_1 \cdots a_i \shS b_1 \cdots b_j }  
		\bkR[a]{ a_{i+1} \cdots a_m \shS b_{j+1} \cdots b_n } 
	}	
	where 
	$a_1,\ldots,a_m,b_1,\ldots,b_n \in \SetO{x,y}$	
	and
	$0 < l \leq m+n$.
This 
	is 
	a spacial case of \refEq{2_PL_DefShS} if $l=1$
	and
	can be proved by induction on $l$.
Let $k=m+n$.
Using 
	\refEq{4_PL_EqShS} with $l=\bkF{\frc[s]{2}{k}}$, 
	we can 
	reduce 
	shuffle products of weight $k$	
	to
	combinations of those of about half weight.
We denote 
	by $\spdO[]{}$ and $\hpdO[]{}$
	created datas of shuffle and stuffle,
	respectively.
They
	map pairs of mult-indices to combinations
	including temporal indeterminates $\bztSO{1,\ldots,1,\bfN}$ $(\bfN \in \midsN \setminus \SetO{\midE} )$
	that are binary versions of regularized MZVs.		
For instance,
	\envHLineCENmPd[4_PL_ExFile1]
	{
	 	\spdO{(1),(2)}
	}
	{
	 	\bztS{1,2} 
	}
	{
		\hpdO{(1),(2)}
	}
	{
		\bztS{1,2} + \bzt{2,1} + \bzt{3}
	}
	{
	 	\spdO{(1,1),(2)}
	}
	{
	 	\bztS{1,1,2} + \bzt{2,1,1}  
	}
	{
	 	\hpdO{(1,1),(2)}
	}
	{
	 	\bztS{1,1,2} + \bztS{1,2,1} + \bzt{2,1,1} + \bztS{1,3} + \bzt{3,1}
	}
In the case of shuffle,
	for using \refEq{4_PL_EqShS} with $l=\bkF{\frc[s]{2}{k}}$,	
	we also create 
	maps 
	from pairs of words to combinations of words
	up to weight $\frc[s]{2}{22} = 11$.
In those additional maps,	
	we allow the words that can not be written in terms of mult-indices (e.g., $x$ and $yx=z_{1}x$).	

Let $\bfK$ be a mult-index 		
	that is expressed as 
	$z_{\bfK}=y^nz_{\bfN}$,
	where
	$n \geq 0$ and $\bfN \in \midsN \setminus \SetO{\midE}$.	
Let $\bfN'$
	denote a mult-index 
	such that $z_{\bfN} = x z_{\bfN'}$.
By \cite[Proposition 8]{IKZ06},	
	\envHLinePd[4_PL_EqHmsR]
	{
		\hmsRO{y^nz_{\bfN}}	
	}
	{
		\mo^n x (y^n  \shS z_{\bfN'}) 
	}	
Since 
	the regularized MZV of  $z_{\bfK}$
	is
	$\evZO[]{} \circ \hmsRO{y^nz_{\bfN}}$,
	its binary version should be	
	\envOTLinePd[4_PL_EqBztR]
	{
		\bztSO{ \bfK }
	}
	{
		\bevZO[]{} \circ \hmsRO{y^n z_{\bfN}}		
	}
	{
		\bevHO[]{} \circ \lnCO{ x (y^n  \shS z_{\bfN'}) } 
	}
The dictionary datas that the file $\fleM[2]$ creates 
	are
	obtained by
	applying \refEq{4_PL_EqBztR}
	to ones that the previous file creates.
For instance,	
	we have by \refEq{4_PL_EqBztR}
	\envHLineCECm
	{
		\bztS{1,2}
	}
	{
		\bevHO[]{} \circ \lnCO{ x (y  \shS y) } 
	\lnP{=}
		\bevHO{ \lnCO{ 2 xyy } }
	\lnP{=}
		0
	}
	{
		\bztS{1,1,2} 
	}
	{
		\bevHO[]{} \circ \lnCO{ x (y^2  \shS y) } 
	\lnP{=}
		\bevHO{ \lnCO{ 3 xyyy } }
	\lnP{=}
		\bzt{2,1,1} 
	}
	{
		\bztS{1,2,1} 
	}
	{
		\bevHO[]{} \circ \lnCO{ x (y  \shS yy) } 
	\lnP{=}
		\bevHO{ \lnCO{ 3 xyyy } }
	\lnP{=}
		\bzt{2,1,1} 
	}
	{
		\bztS{1,3}
	}
	{
		\bevHO[]{} \circ \lnCO{ x (y  \shS xy) } 
	\lnP{=}
		\bevHO{ \lnCO{ xyxy + 2xxyy } }
	\lnP{=}
		\bzt{2,2} 
	}
	and so 
	the previous datas in \refEq{4_PL_ExFile1}
	are converted to 
	\envHLineCECmNm[4_PL_ExFile2]
	{
	 	\spdRO{(1),(2)}
	}
	{
	 	0
	}
	{
		\hpdRO{(1),(2)}
	}
	{
		\bzt{2,1} + \bzt{3}
	}
	{
	 	\spdRO{(1,1),(2)}
	}
	{
	 	0
	}
	{
	 	\hpdRO{(1,1),(2)}
	}
	{
	 	\bzt{2,1,1} + \bzt{2,2} + \bzt{3,1}
	}
	where
	$\spdRO[]{}$ and $\hpdRO[]{}$ 
	stand for  
	the maps 	
	of regularized shuffle and stuffle,	
	respectively.

The file $\fleM[3]$ 
	makes binary EDS relations,	
	\envHLineCm
	{
		\hpdRO{\bfK,\bfL} + \spdRO{\bfK,\bfL}
	}
	{
		0
		\qquad
		( (\bfK,\bfL) \in \epidsO{k} )	
	}
	by combining 
	the previous datas. 		
For instance,
	the datas in \refEq{4_PL_ExFile2}
	create the two relations,
	\envHLineCm		
	{
		\bzt{2,1} + \bzt{3}
	}
	{
	 	0
	}
	\vPack[25]\envHLinePd
	{
	 	\bzt{2,1,1} + \bzt{2,2} + \bzt{3,1}
	}
	{
		0
	}

The file $\fleM[4]$ 
	converts 
	the binary EDS relations of a weight $k$
	to
	a binary linear system (which we call a binary EDS linear system)
	in both of text and pickle formats.
The text format
	is organized as follows:
\envItemEm{
\item
A line with the first character `\#' 
	is a comment line.
Comment lines typically occur at the beginning of the file, 
	but are allowed to appear throughout the file.
\item
The remainder of the file 
	contains 
	lines defining the binary linear relations, 
	one by one.
\item
A relation is 
	defined by
	positive integers numbering binary MZSs.
A number `0' is 
	typically 
	placed at the last of the line,
	but 
	it is optional.
}
For example,
	the line ``2 4 0''
	is corresponding to $\bzt{2,1}+\bzt{3}=0$,
	if $\bzt{2,1}$ and $\bzt{3}$ 
	are 
	numbered as $2$ and $4$,
	respectively.
The pickle files are not necessary but useful:
	e.g., when loading the large size system.	
For \refExp{3_EXP3},
	we also make 
	binary $\letKNT$ and $\letMJPO$ linear systems
	by restricting 
	binary EDS relations. 	

The above programs run under the parallel process 
	since
	the datas can be created independently
	if	
	$\epidsO{k}$ is divided into a plurality of blocks.
The filenames of the datas by the parallel process
	have 
	strings `$\_\mathtt{Bn}$' $(\mathtt{n} \in \setN)$ at their tails.	
Editing the file $\fleConfig$
	we can control the max number of parallel threads. 

In \refTab{4_Tbl_Main_make_tme},
	we present 
	computation times (or elapsed real times) 
	to execute all files for $k\geq18$,
	where 
	$\FLEm[i]$ 
	stands for the $i$-th file mentioned above from $0$ to $4$.
We find that
	calculating the regularizations in $\FLEm[2]$
	is 
	the dominant process.
\refTab{4_Tbl_Main_make_file}
	lists
	the file sizes of the binary linear systems in pickle format
	for $\epidsKO{k}$, $\epidsMO{k}$ and $\epidsO{k}=\epidsEO{k}$.
As expected from 
	$\epidsKO{k} \subset \epidsMO{k} \subset \epidsO{k}$,
	the file of $\letKNT$ is smallest 
	and
	that of EDS is largest for each weight.
The size of text format file 
	is 
	about $1.5$ times the size of pickle one.
For each weight $k$,
	the maximum memory size (or maximum resident set size) 
	to execute the files $\FLEm[i]$ 
	is
	the size required by $\FLEm[4]$, 	
	which is
	about half size used in Gaussian forward elimination 	%
	(see \refTab{4_Tbl_Main_cal_bspZ}).	
Computationally,	
	making linear systems is not harder than calculating their coranks (or dimensions of cokernel) 
	as we will see below.		

\begin{table}[!t]\renewcommand{\arraystretch}{1.3} \tabcolsep=10pt	
\begin{center}
	\caption{
		Elapsed real times [sec] to make binary systems.
	}\label{4_Tbl_Main_make_tme}
	\vWiden
	\begin{tabular}{|c|ccccc|c|}\hline 
		$k$ & $\FLEm[0]$ & $\FLEm[1]$ & $\FLEm[2]$ & $\FLEm[3]$ & $\FLEm[4]$ & Total \\\hline 
		18 & 0 & 75 & 246 & 53 & 54 & 428 ($\fallingdotseq$ 7min) \\\hline 
		19 & 1 & 188 & 805 & 133 & 186 & 1313 ($\fallingdotseq$ 22min) \\\hline 
		20 & 3 & 469 & 3137 & 510 & 543 & 4662 ($\fallingdotseq$ 1.3hour) \\\hline 
		21 & 7 & 1529 & 15607 & 1384 & 2362 & 20889 ($\fallingdotseq$ 5.8hour) \\\hline 
		22 & 15 & 3018 & 61898 & 3675 & 6578 & 75184 ($\fallingdotseq$ 21hour)\\\hline 
	\end{tabular}
\end{center}
\end{table}
\begin{table}[!t]\renewcommand{\arraystretch}{1.3} \tabcolsep=10pt	
\begin{center}
	\caption{
		File sizes of binary linear systems in pickle format.
	}\label{4_Tbl_Main_make_file}
	\vWiden
	\begin{tabular}{|c|c|c|c|}\hline 
		$k$ & $\letKNT$ &  $\letMJPO$ & EDS \\\hline 
		18 & 8.3M & 150M & 274M \\\hline 
		19 & 21M & 509M & 922M \\\hline 
		20 & 47M & 1.6G & 2.8G \\\hline 
		21 & 105M & 5G & 8.6G \\\hline 
		22 & 233M & 16G & 26G \\\hline 
	\end{tabular}
\end{center}
\end{table}

\subsection{Executable file in $\dirC$}
The file $\fleC$ ($\FLEc[d]$ for short)
	executes the Gaussian forward elimination
	on a given binary linear system of a weight $k$ 
	by using \refAlg{A_ALG1}.
In the process,
	an order of mult-indices (or binary MZSs)
	have to be determined 
	to 
	convert 
	the inputted binary linear system into the corresponding binary matrix.	
We employ 	
	a sequence $({\bfK_1}, \ldots, {\bfK_{2^{k-2}}})$		
	satisfying the following:	
	if $i<j$,	
\envItem{
\item[\cndBO{1}]
	$\dpt{\bfK_i} > \dpt{\bfK_j}$; or
\item[\cndBO{2}]
	$\dpt{\bfK_i} = \dpt{\bfK_j}$
	and
	$(\bfK_i,\bfK_j) \notin \midsHO{k} \times (\midsO{k}\setminus\midsHO{k})$.
}
The condition \cndBO{1}
	means that
	the mult-indices (or columns in the corresponding matrix) are sectioned into $k-1$ blocks by depth:
	the mult-indices in a left block have 
	a greater depth than those in a right block.
The condition \cndBO{2}
	means that
	the Hoffman mult-indices of depth $r$
	are
	at the rightmost place in the $(k-r)$th block. 	
For example,
	the order of weight $4$ 
	determined by
	$\bfK_1=(2,1,1)$, $\bfK_2=(3,1)$, $\bfK_3=(2,2)$ and $\bfK_4=(4)$
	satisfies \cndBO{1} and \cndBO{2}:
	they are sectioned as $(\bfK_1) | (\bfK_2,\bfK_3) | (\bfK_4)$
	and
	the only Hoffman mult-index $\bfK_3$ 
	is located rightmost in the $2$th block.

\refPrp{4_PRP1} is shown as follows.	

\envProof[\refPrp{4_PRP1}]{	
We consider the situation where 
	we run $\FLEc[d]$ 	
	by inputing 
	the binary EDS linear system of weight $k$.	
We then obtain 	
	a row echelon matrix satisfies the following.		
\envItem{
\item[\cndEO{1}]
	There exists a non-zero pivot at any column $\bfK$ in $\midsO{k}\setminus\midsHO{k}$.
\item[\cndEO{2}]
	There exists no non-zero pivot at any column $\bfK$ in $\midsHO{k}$.
}
%
For a non-zero combination $c = \bZtO{\bfK_{i_1}} + \cdots + \bZtO{\bfK_{i_j}}$ 
	in $\bspUO{k}$
	with $i_1 < \cdots < i_j$,
	we define the leading term of $c$ 	
	by 
	\envHLinePd
	{
		\ltm{c}
	}
	{
		\bZtO{\bfK_{i_1}}
	}
By \cndBO{1} and \cndBO{2},	
	the statements \cndEO{1} and \cndEO{2}
	are equivalent to 
	\cndEdO{1} and \cndEdO{2},
	respectively:
\envItem{
\item[\cndEdO{1}]
	There exists 
	a combination $c \in \bspEO{k}$
	such that
	$\ltm{c} = \bZtO{\bfK}$ 	
	for any $\bfK$ in $\midsO{k}\setminus\midsHO{k}$.
\item[\cndEdO{2}]
	There exists 
	no combination $c \in \bspEO{k}$
	such that
	$\ltm{c} = \bZtO{\bfK}$ 
	for any $\bfK$ in $\midsHO{k}$.
}
Under \cndEdO{1} and \cndEdO{2},	
	the back substitution (performed imaginarily) 	
	implies 
	\refPrp{4_PRP1},
	where
	the fact that
	$\midsHT{k}{r} = \phi $ for any depth $r < \bkF{\frc[s]{3}{k}}$
	is 
	used for \refEq{4_PRP1_IncBHS}.	
}

We give
	examples of \refEq{4_PL_IncBZSq} for $k\leq7$
	excluding the case that $\bztqO{\bfK} = 0$.
Note that
	$\bztqO{\bfK}$ is always zero	
	if
	$\bfK \in \midsT{k}{r}$
	and
	$\midsHT{k}{r} = \phi$. 
	\envPLine{\HLineCECm[p]
	{
		\bztqO{3,1}
	}
	{
		\bztqO{2,2}
	}
	{\rule{0pt}{25pt}
		\bztqO{4,1}
	}
	{
		\bztqO{2,3} + \bztqO{3,2} 
	}
	{\rule{0pt}{25pt}
		\bztqO{2,1,3}	\lnP{=}	\bztqO{3,2,1}	\lnP{=}	\bztqO{4,1,1}
	}
	{
		\bztqO{2,2,2}
	}
	{
		\bztqO{5,1}	
	}
	{
		\bztqO{3,3}
	}\HLineCTCm[p]
	{\rule{0pt}{25pt}
		\bztqO{5,1,1}	\lnP{=}	\bztqO{3,1,3}	
	}
	{
		\bztqO{2,2,3} + \bztqO{2,3,2} + \bztqO{3,2,2}
	}\HLineCEPd
	{
		\bztqO{3,3,1}
	}
	{
		\bztqO{2,3,2}
	}
	{
		\bztqO{4,2,1}
	}
	{
		\bztqO{2,2,3} + \bztqO{2,3,2}
	}
	{
		\bztqO{4,1,2}	\lnP{=}	\bztqO{2,1,4}
	}
	{
		\bztqO{2,2,3} + \bztqO{3,2,2}
	}
	{
		\bztqO{2,4,1}
	}
	{
		\bztqO{2,3,2} + \bztqO{3,2,2}
	}
	}
Examining the Gaussian forward elimination performed by $\FLEc[d]$ in detail,
	we can find 
	a part of the inputted binary EDS relations which forms a basis of $\bspEO{k}$.
We give examples of bases for $k\leq6$,
	where
	only the pairs of mult-indices are written
	(see \refTab{2_Tbl_EDS} that lists associated relations for $k\leq4$).
	
\envItem{ \setlength{\leftskip}{10pt}
\item[\ulTx{$k=3$}]
	$((1), (2))$.
\item[\ulTx{$k=4$}]
	$((1),(2,1))$, $((1),(3))$, $((2),(2))$.
\item[\ulTx{$k=5$}]
	$((1),(2,1,1))$, $((1),(2,2))$, $((1),(3,1))$, $((1),(4))$, $((2),(2,1))$, $((2),(3))$.
\item[\ulTx{$k=6$}]
	$((1),(2,1,1,1))$, $((1),(2,1,2))$, $((1),(2,2,1))$, $((1),(2,3))$, $((1),(3,1,1))$, $((1),(3,2))$, 
	$((1),(4,1))$, $((1),(5))$, $((2),(2,1,1))$, $((2),(2,2))$, $((2),(3,1))$, $((2,1),(2,1))$, 
	$((2,1),(3))$, $((3),(3))$.
}

We can verify 
	\refExp{3_EXP3} similarly to \refExp{3_EXP2}.
We input 
	$\letKNT$ and $\letMJPO$ linear systems 
	into $\FLEc[d]$.
By \refTab{3_Tbl_EXP3},
	in most cases,
	row echelon matrices that do not satisfy \cndEO{1} 
	are outputted.
The fails of \cndEO{1} 	
	induce \refEq{3_EXP3_Eq},
	and 
	ensure \refExp{3_EXP3}.		

The program in $\FLEc[d]$ applies the parallel process 
	to determine an order of mult-indices 
	since 
	mult-indices can be divided by depth.
For instance,
	$(k-1)$ parallel threads occur as preprocessing 
	if a binary EDS linear system of weight $k$ is inputted.
\refAlg{A_ALG1}, 
	the main process for computing a row echelon matrix, 
	is executed in single.
It appears that
	the parallelization of \refAlg{A_ALG1}  is not easy 
	because 
	a non-simple search procedure is incorporated.	

In \refTab{4_Tbl_Main_cal_bspZ},
	we present 
	the statics of the executions by $\FLEc[d]$	
	whose inputs are 
	the binary $\letKNT$, $\letMJPO$ and EDS relations.
We observe that
	the computation for $\letKNT$ 	
	requires 
	much more time than $\letMJPO$ and EDS,
	although
	the number of relations of $\letKNT$ is 
	quite small
	such that the corresponding matrix is square for any $k \geq 7$.
This phenomenon expresses 
	a characteristic of \refAlg{A_ALG1}. 	
It employs   
	a conflict based search procedure	
	inspired by 
	the conflict-driven clause learning (CDCL),
	a modern method with many successes to practical applications
	in solving the Boolean satisfiability (SAT) problem.	
Roughly speaking,
	relations with good structures for finding conflict combinations		
	can accelerate searching a pivot relation	
	(see \refRem{A_Rem1} for more information).
The memory cost is bad in comparison with the statics in \cite{KNT08},
	but 
	the runtime is about 10 times more faster.
Therefore	
	we can improve the record of calculating \refEq{3_PL_CnjEDS} from $k=20$ to $22$
	by the use of
	a machine with large memory capacity.	

\begin{table}[!t]\renewcommand{\arraystretch}{1.3} \tabcolsep=10pt	
\begin{center}
	\caption{
		Statistics of the computations of \refExp[s]{3_EXP2} and \ref{3_EXP3}.
		`Rels' is the number of relations.
		`MeanNum' is the average number of terms per relation. 
		`Memory' and `Time' are
		the resident set size and elapsed real time,
		respectively.
		In each block with respect to the weight $k$,
		top row indicates information on $\letKNT$,
		middle row indicates that on $\letMJPO$ 
		and
		bottom row indicates that on EDS. 
	}\label{4_Tbl_Main_cal_bspZ}
	\vWiden	
	\begin{tabular}{|c|ccccc|l}\cline{1-6} 
		$k$ & $2^{k-2}$ & Rels & MeanNum & Memory & Time & \\\cline{1-6} 		
		\multirow{3}{*}{18}  & \multirow{3}{*}{65536} & 65536 & 30.1 & 4.6G & 8.6hour & $\letKNT$\\\cdashline{3-7}  		
		  &  & 155711 & 230.4 & 7.3G & 8.8min & $\letMJPO$  \\\cdashline{3-7} 
		 &  & 188470 & 364.4 & 11.4G & 9.8min & EDS  \\\cline{1-6}		  
		\multirow{3}{*}{19}  & \multirow{3}{*}{131072} & 131072 & 33.7 & 16.5G & 68hour & \\\cdashline{3-6} 			  
		  &  & 327679 & 339.5 & 22.9G & 42.4min &  \\\cdashline{3-6} 
		  &  & 393206 & 523.1 & 34.3G & 43.7min &  \\\cline{1-6} 		
		\multirow{3}{*}{20}  & \multirow{3}{*}{262144}  & 262144 & 37.6 & 61G & 22day & \\\cdashline{3-6}  					
		  &  & 688254 & 500.5 & 82G & 5.3hour &  \\\cdashline{3-6} 		
		 & & 819316 & 751.7 & 110G & 4.7hour &  \\\cline{1-6}
		\multirow{3}{*}{21}  & \multirow{3}{*}{524288} & - & - & - & - & \\\cdashline{3-6}
		  &  & 1441791 & 739.8 & 256G & 30hour &  \\\cdashline{3-6} 			
		 &  & 1703925 & 1083.3 & 329G & 25hour &  \\\cline{1-6} 	 
		\multirow{3}{*}{22}  & \multirow{3}{*}{1048576} & - & - & - & - & \\\cdashline{3-6}  		
		  &  & 3014911 & 1094.4 & 789G & 8day &  \\\cdashline{3-6} 	
		 &  & 3539188 & 1564.1 & 982G & 7day &  \\\cline{1-6}		
	\end{tabular}
\end{center}
\end{table}

\section{Problem} \label{sectFive}
Some problems arise 
	in connection with the experiments in \refSect{sectThree}.

\refExp[s]{3_EXP1} and \ref{3_EXP2}
	indicate  
	typical problems on the dimensions of $\bspZO{k}$ and $\bspZqT{k}{r}$:	
	obviously,
	\refPrb{5_PRB2} includes \refPrb{5_PRB1}.

\begin{problem}\label{5_PRB1}
Does \refEq{3_EXP1_Eq} hold for any weight $k$?
\end{problem}
\begin{problem}\label{5_PRB2}
Does \refEq{3_EXP2_Eq} 
	(or \refPrp{4_PRP1}) hold for any weight $k$ and depth $r$?
\end{problem}

\refExp{3_EXP3}		%
	yields
	the following:	
\begin{problem}\label{5_PRB3}
{\bf (i)}
Is there a subset $M \subset \setN$ 	
	such that 
	$ M \cap \bkS{22} = \SetO{7,15} $		
	and		
	the sequence $\bkR{\cdmBZMO{k}}_{k\geq0}$ satisfies \refEq{3_PL_EqFibMJPO}?
\mbox{}\\{\bf (ii)}
Can we find a law in the sequence $\bkR{\cdmBZKO{k}}_{k\geq0}$?
\end{problem}

We have adopted the binary field $\setF[2]$ 
	for the scalar field of the formal multiple zeta space
	and
	for the computation of corank.
(It is worth noting that
	the experiments of \cite{KNT08} 
	employ  
	$\setF[16381]$ and $\setF[31991]$.)
There are no particular reasons for choosing $\setF[2]$
	except 
	computational science techniques 	
	are easy to apply.
A discovery of a regularity of $\cdmBZqT{k}{r}$ in \refTab{3_Tbl_EXP12_dd}		
	is 	
	a product of good luck.	

\begin{problem}\label{5_PRB4}
{\bf (i)}
Can we find a theoretical reason 
	why 
	the dimensions $\cdmBZqT{k}{r}$ $(r<k\leq 22)$ 
	have a Pascal triangle pattern?
\mbox{}\\{\bf (ii)}
What will the dimensions be 
	if
	we adopt other finite fields $\setF[p]$?
\end{problem}

Like MZVs,
	we can make an assumption that
	binary MZSs satisfy a multiplication compatible 
	with the shuffle and stuffle products.	
Under the assumption,
	we have
	\envM
	{
		\bztO{2}^2  
	}
	{
		0
	}
	since		
	\envMFPd
	{
		\bevZO{z_2} \bevZO{z_2}
	}
	{
		\bevZO{ z_2 \shS[\,] z_2 }
	}
	{
		\bevHO[]{} \circ \lnCO{ 2 z_{2,2} + 4 z_{3,1} }
	}
	{
		0
	}	
This means that
	the algebras of MZV and binary MZS 
	are 
	different.
In particular,
	$\setF[2] \bkS{\bztO{2}} = \SpO{1, \bztO{2}}_{\setF[2]}$
	is 
	not isomorphic to the polynomial ring in one variable,
	and
	statements and conjectures involving $\spZN/\zt{2}\spZN$ 
	(e.g., those involving finite and symmetric multiple zeta values introduced in \cite{Kaneko19})
	can not be varied 
	to $\bspZN/\bzt{2}\bspZN$ directly.
It seems a mysterious problem
	that
	whether the algebra of binary MZS has 
	a good property 
	and
	a connection to the algebra of MZV.


\section*{Acknowledgements}
The author would like to 
	thank 
	Tomohiro Sonobe
	for help with computing environments,
	and
	Junichi Teruyama 
	for a recommendation to use \refEq{4_PL_EqShS}
	which 
	made it possible to reduce computation costs. 
This work was supported by 
	Japan Society for the Promotion of Science,
	Grant-in-Aid for Scientific Research (C) 20K03727. 
	

\section*{Appendix}
We will introduce 
	a technique to speed up the Gaussian forward elimination over any field $K$
	for 
	a system of linear combinations that have some structure.		
An essential part of the technique 
	appears in \cite{MS18}
	to decide the full rankness of a binary matrix.

Let $x_1, \ldots, x_n$ be variables,
	and
	we order the variables according to their subscripts.		
For a non-zero linear combination 
	$p=\Fc{p}{x_1,\ldots,x_n}=\SmO{}{} c_i x_i$ over $K$,
	we 
	denote by $\sbsFO{l}$ and $\coeFO{l}$
	the subscript and coefficient of the minimum variable,	
	respectively.
That is,
	$\sbsFO{p}=\MinT{i}{ c_i \neq0}$
	and
	$\coeFO{p}=c_{\sbsFO{p}}$.
We define
	$\sbsFO{p}=n+1$ and $\coeFO{p}=0$
	when
	$p=0$.

In what follows
	we will handle 
	mainly linear combinations over $K$,
	and 
	we just call them combinations.
Let $\cmbsO[]{p_1,\ldots,p_m}$
	denote the $K$-vector space spanned by combinations $p_1,\ldots,p_m$,
	and
	let 
	$\cmbsO{p_1,\ldots,p_m} = \cmbsO[]{p_1,\ldots,p_m} \setminus \SetO{0}$.
We say that
	$p_i$ is a pivot combination if $\sbsFO{p_{i}} = i$,
	and	
	\envHLineT
	{
		(p_{i_g})_{1 \leq g \leq h}
	}
	{
		(p_{i_1},\ldots,p_{i_h})
	} 	
	is
	a pivot sequence 
	if
	$1 \leq i_1 < \cdots < i_h \leq n$	
	and
	every $p_{i_g}$ is a pivot combination.	

There are two key processes for the speed-up technique.
One is 
	a conflict search procedure.	
\begin{processA}\label{A_PRC1}\mbox{}\vPack[5]
\envItemDp{\setlength{\itemsep}{0pt}		
\item[]Input:	
	Combinations $L = \SetO{l_1,\ldots,l_m}$
	and
	a pivot sequence $(p_1,\ldots,p_{j-1})$.	
\item[]Output:
	Either $(0,\midE)$
	or
	$(q_i,\bfK_i)$ 
	such that
	\vPack[5]\envItemEm{\setlength{\itemsep}{0pt}		
	\item[(a)]
		$q_i\in L$	
		with
		$\sbsFO{q_i}=i\leq j$;
	\item[(b)]
		$\bfK_i = (k_i,\ldots,k_{j-1},k_j,\ldots,k_n) \in K^{n-i+1}$ with $k_j=1$ and $k_{j+1}= \cdots = k_n = 0$; 
	\item[(c)]
		$q_i \mVert_{(x_i,\ldots,x_n)=\bfK_i} \in K^*$; and		
	\item[(d)]
		$p_i \mVert_{(x_i,\ldots,x_n)=\bfK_i} = \cdots = p_{j-1} \mVert_{(x_{i},\ldots,x_n)=\bfK_i} = 0$.		
	}
}
\envItemEm{\renewcommand{\labelenumi}{\arabic{enumi}.}	
\item
	Set $\bfK_j=(k_j,\ldots,k_n)=(1,0,\ldots,0) \in K^{n-j+1}$ and $i=j$.
\item
	Search $q_i$ from $\SetT{l \in L}{\sbsFO{l} = i}$ such that $q_i \mVert_{(x_i,\ldots,x_n)=\bfK_i} \in K^*$.
\item
	Return $(q_i,\bfK_i)$ if such $q_i$ exists.
\item
	Return  $(0,\midE)$ if $i=1$.
\item 
	Evaluate $k_{i-1} = - \frc[d]{\coeFO{p_{i-1}}}{ p_{i-1} - \coeFO{p_{i-1}} x_{i-1}} \mVert[g]_{(x_i,\ldots,x_n)=\bfK_i}\in K$.\footnote{
This evaluation 
	is well-defined 
	since $\sbsFO{p_{i-1}}=i-1$ and $\coeFO{p_{i-1}}\neq0$.
The condition (d) 
	follows from 
	\envOTLinePd
	{	
		p_{i-1} \mVert_{(x_{i-1},\ldots,x_n)=(k_{i-1},\ldots,k_n)} 
	}{
		\coeFO{p_{i-1}} k_{i-1} 
		+ 
		\bkR{ p_{i-1} - \coeFO{p_{i-1}} x_{i-1} } \mVert_{(x_{i},\ldots,x_n)=(k_{i},\ldots,k_n)}
	}{
		0
	}
}
\item
	Set $\bfK_{i-1} = (k_{i-1},\bfK_i)$.
\item
	Update $i \leftarrow i-1$,
	and
	go back to step $2$.
}
\end{processA}		

Another 
	is 
	the classical elimination procedure with an evidence of conflict.
\begin{processA}\label{A_PRC2}\mbox{}\vPack[5]
\envItemDp{\setlength{\itemsep}{0pt}		
\item[]Input:	
	A pivot sequence $(p_i,\ldots,p_{j-1})$
	and
	a pair $(q_i,\bfK_i)\neq (0,\midE)$ 
	which satisfies the output conditions in \refPrc{A_PRC1}.
\item[]Output:
	A combination $q_j\in\cmbsO{q_i,p_i,\ldots,p_j}$ 
	such that $\sbsFO{q_j}=j$.\footnote{
The theory of Gaussian elimination 
	only ensures 
	$q_j\in\cmbsO[]{q_i,p_i,\ldots,p_j}$ and $\sbsFO{q_j} \geq j$.
However,
	updating method of $q$ in step 2,	
	together with the output conditions (c) and (d) in \refPrc{A_PRC1},
	implies 
	\envMPdPt{\in}
	{
		q_j \mVert_{(x_i,\ldots,x_n)=\bfK_i} 
	}
	{
		K^*
	}	
It also implies 	
	$\sbsFO{q_j} = j$.
In fact,
	if $\sbsFO{q_j} > j$,
	\envOFLineCm
	{
		q_j \mVert_{(x_i,\ldots,x_n)=\bfK_i} 
	} 
	{
		q_j \mVert_{(x_{j+1}, \ldots, x_n)=(k_{j+1}, \ldots, k_n)}
	}
	{
		q_j \mVert_{x_{j+1}= \cdots = x_n=0}
	}
	{
		0
	}
	which
	is a contradiction.	
Therefore 	 
	the output condition in \refPrc{A_PRC2}
	holds.	
}
}
\envItemEm{\renewcommand{\labelenumi}{\arabic{enumi}.}	
\item
Set $q=q_i$.
\item 
For $h$ from $i$ to $j-1$,
	update $q \leftarrow q - \frc[d]{\coeFO{p_h}}{\coeFO{q}}p_h$	
	if $h = \sbsFO{q}$.
\item
Return $q_j=q$.
}
\end{processA}		

We can construct
	a process
	to find a new pivot combination
	by
	combining 
	\refPrc[s]{A_PRC1} and \ref{A_PRC2}.	

\begin{processA}\label{A_PRC3}\mbox{}\vPack[5]	
\envItemDp{\setlength{\itemsep}{0pt}		
\item[]Input:	
	Combinations $l_1,\ldots,l_m$
	and
	a pivot sequence $(p_1,\ldots,p_{j-1})$.	
\item[]Output:
	Either 
	$0$ 
	or 
	a combination $p_j\in\cmbsO{l_1,\ldots,l_m,p_1,\ldots,p_j}$ 
	such that $\sbsFO{p_j}=j$.		
}
\envItemEm{\renewcommand{\labelenumi}{\arabic{enumi}.}	
\item
Receive $(q_i,\bfK_i)$ 
	from \refPrc{A_PRC1} 
	for the inputs $L=\SetO{l_1,\ldots,l_m}$ and $(p_1,\ldots,p_{j-1})$.
\item
Return $0$ if $q_i=0$.
\item
Receive $q_j$ from \refPrc{A_PRC2}
	for the inputs  $(p_i,\ldots,p_{j-1})$ and $(q_i,\bfK_i)$.
\item
Return $p_j=q_j$.
}
\end{processA}		

\refPrc{A_PRC3}
	is 
	essential 
	for
	finding a pivot combination whose minimum variable is $x_j$,
	because
	we can find out it by \refPrc{A_PRC3}	
	if and only if 
	it exists. 
	
\begin{propositionA}\label{A_PRP1}
For combinations $l_1, \ldots, l_m$ 
	and
	a pivot sequence $(p_1,\ldots,p_{j-1})$,
	the following statements are equivalent.
\mbox{}\\{\bf(i)}
\refPrc{A_PRC3}
	outputs 
	$p_j \in \cmbsO{l_1,\ldots,l_m,p_1,\ldots,p_{j-1}}$ such that $\sbsFO{p_j}=j$.
\mbox{}\\{\bf(ii)}
There exists a combination  
	$p_j \in \cmbsO{l_1,\ldots,l_m,p_1,\ldots,p_{j-1}}$ such that $\sbsFO{p_j}=j$.	
\end{propositionA}	
\envProof{
Obviously,
	{\bf(i)} implies {\bf(ii)}.
Suppose 
	{\bf(ii)} is true to prove the converse.
Then there exist elements 
	$c_1,\ldots,c_m, d_1,\ldots, d_{j-1}$ in $K$	
	such that
	\envHLinePd
	{
		p_j
	}
	{
		\SmO{h} c_h l_h
		+
		\SmO{i} d_i p_i
	}
We have
	$p_j\mVert_{(x_j,x_{j+1},\ldots,x_n)=(1,0,\ldots,0)} \in K^*$
	since $x_j$ is the minimum variable in $p_j$.	
	
We first consider the situation where		
	we run \refPrc{A_PRC1}
	for the inputs $L=\SetO{l_1,\ldots,l_m}$ and $(p_1,\ldots,p_{j-1})$:
	however,		
	we temporally assume that
	step 3 is skipped  
	and 
	the process ends with the output $(0,\midE)$ at step 4 of $i=1$.	
Let $k_1,\ldots,k_{j-1}$ 
	be 
	the elements in $K$
	which are 
	recursively determined   
	as at step $5$,
	and
	let 
	$\bfK = (k_1,\ldots,k_{j-1},1,0,\ldots,0) \in K^n$.
Then		
	$\Fc{p_1}{\bfK}=\cdots=\Fc{p_{j-1}}{\bfK}=0$,
	and
	\envOTLinePd
	{
		\Fc{p_j}{\bfK}
	}
	{
		\SmO{h} c_h \Fc{l_h}{\bfK} + \SmO{i} d_i \Fc{p_i}{\bfK}
	}	
	{
		\SmO{h} c_h \Fc{l_h}{\bfK}
	}
Since
	$\Fc{p_j}{\bfK} = p_j\mVert_{(x_j,x_{j+1},\ldots,x_n)=(1,0,\ldots,0)} \in K^*$,
	this	
	implies 
	$\Fc{l_h}{\bfK} \in K^*$ for some $h$,
	which means that	
	\refPrc{A_PRC1} 
	can 
	find out $q_i$ in step 2 such that $q_i \mVert_{(x_i,\ldots,x_n)=\bfK_i} \in K^*$,
	at least when $i=\sbsFO{l_h}$.
Therefore,
	\refPrc{A_PRC1} without the temporal assumption 	
	always outputs 
	$(q_i,\bfK_i) \neq (0,\midE)$. 

We input $L=\SetO{l_1,\ldots,l_m}$ and $(p_1,\ldots,p_{j-1})$ into \refPrc{A_PRC3}.		
At step $1$,
	we receive $(q_i,\bfK_i)\neq (0,\midE)$
	from \refPrc{A_PRC1}.	
Thus
	step 2 is skipped,
	and
	$q_j$ is received from \refPrc{A_PRC2} at step 3,
	which satisfies
	the condition required in {\bf(i)}.
Since
	$p_j=q_j$ is returned at step 4,
	we conclude  
	{\bf(i)} holds.	
}

For a subscript $j$ and a pivot sequence $(p_{i_g})=(p_{i_g})_{1 \leq g \leq h}$ with $i_h < j$, 
	we define 
	\envHLineDefPd
	{
		\sbssDT{(p_{i_g})}{j}
	}
	{ 
		\bkS{j-1}	\setminus \SetO{i_1,\ldots,i_h}
	}	
We call an integer in $\sbssDT{(p_{i_g})}{j}$ a deficient subscript,
	and
	a variable $x_i$ with $i\in\sbssDT{(p_{i_g})}{j}$ a deficient variable.
	
We need to modify \refPrc{A_PRC3}
	for practical use.
	
\begin{processA}\label{A_PRC4}\mbox{}\vPack[5]
\envItemDp{\setlength{\itemsep}{0pt}		
\item[]Input:	
	Combinations $l_1, \ldots, l_m$,
	a subscript $j$,
	and
	a pivot sequence $(p_{i_g})_{1 \leq g \leq h}$ with $i_h < j$.	
\item[]Output:
	Either 
	$0$
	or
	a combination $p_j\in\cmbsO{l_1,\ldots,l_m,p_{i_1},\ldots,p_{i_h}}$ 
	with
	$\sbsFO{p_j} \in \sbssDT{(p_{i_g})}{j} \cup \SetO{j}$.
}
\envItemEm{\renewcommand{\labelenumi}{\arabic{enumi}.}	
\item 
	Change the variable order  
	by 
	moving 	
	the deficient variables backward.
\item
Prepare 
	the pivot sequence $(p'_1,\ldots,p'_{j-1-\nmSet{\sbssDT{(p_{i_g})}{j} }})$
	for
	the new variable order.  
\item
Receive $p'_{j-\nmSet{\sbssDT{(p_{i_g})}{j}}}$ from \refPrc{A_PRC3} 
	for the inputs $l_1,\ldots,l_m$ and $(p'_1,\ldots,p'_{j-1-\nmSet{\sbssDT{(p_{i_g})}{j}}})$.
\item
Undo the variable order 
	by putting the deficient variables back to their original places.
\item
Return $p_j=p'_{j-\nmSet{\sbssDT{(p_{i_g})}{j}}}$.
}
\end{processA}		

We are in a position to state \refAlg{A_ALG1}	
	for a fast Gaussian forward elimination.

\begin{algorithmA}\label{A_ALG1}\mbox{}\vPack[5]
\envItemDp{\setlength{\itemsep}{0pt}		
\item[]Input:	
	Combinations $L=\SetO{l_1, \ldots, l_m}$. 
\item[]Output:
	A pivot sequence $(p_{i_g})$. 	
}
\envItemEm{\renewcommand{\labelenumi}{\arabic{enumi}.}	
\item		
	Create subsets 
	$L_i = \SetT{l\in L}{ \sbsFO{l} = i}$ $(i=1,\ldots,n)$.
\item
Set $j=0$ and $(p_{i_g})=\phi$.		
\item
Execute the following loop process 
	to make a pivot sequence $(p_{i_g})$:
\envItemEm{\renewcommand{\labelenumii}{(\roman{enumii})}
\item[(i)]
Update $j \leftarrow j+1$ if $j < n$;
	otherwise break.	
\item[(ii)]
If $L_j \neq \phi$,
	append a combination in $L_i$ to $(p_{i_g})$
	and
	go back to (i).
\item[(iii)]
Receive $p_j$ from \refPrc{A_PRC4} 
	for the inputs $L_1 \cup \cdots \cup L_{j-1}$ and $(p_{i_g})$.
\item[(iv)]
If $p_j=0$,
	go back to (i).
\item[(v)]
Append $p_j$ to $(p_{i_g})$,\footnote{
The loop process 
	ensures $\sbsFO{p_j}=j$.		
To show this,
	we may prove $\sbsFO{p_j} \notin \sbssDT{(p_{i_g})}{j}$ 
	by the output condition in \refPrc{A_PRC4}.
Suppose $\sbsFO{p_j}  \in \sbssDT{(p_{i_g})}{j}$
	and
	set $j'=\sbsFO{p_j} < j$.
Then,
	on the $j'$-round in the loop process,
	\refPrc{A_PRC4} at (iii)	
	must return 
	a non-zero combination 
	by \refPrp{A_PRP1} and the existence of $p_j$,
	where 
	note that
	\refPrc{A_PRC4} is essentially \refPrc{A_PRC3}.
This means
	a combination $p_{j'}$ satisfying $\sbsFO{p_{j'}} =j'$ 
	must be appended to $(p_{i_g})$ at (v) on the $j'$-round,
	which 
	contradicts $j'=\sbsFO{p_j} \in \sbssDT{(p_{i_g})}{j}$.
}
	and 
	back to (i).
}
\item
Return $(p_{i_g})$.
}
\end{algorithmA}	

The pivot sequence $(p_{i_g})$ outputted 	by \refAlg{A_ALG1}
	is 
	a row echelon matrix under the order $x_1 < \cdots < x_n$
	thanks to \refPrp{A_PRP1}
	(see the footnote in (v) of step 3 for details).

\begin{remarkA}\label{A_Rem1}
\refPrc{A_PRC1}
	is influenced by the unit propagation (UP)
	in
	the algorithm to solve the Boolean satisfiability (SAT) problem
	(see, e.g., \cite[Chapter 1]{BHMW09}).	
SAT is the first problem that was proved to be NP-complete,
	which means that 
	all NP-problems are at most as difficult as SAT.	
UP is 
	a technique 
	to 
	determine an assignment value for the variable we watch
	while 
	searching a conflict combination (or a conflict clause in SAT terminology).

\refPrc{A_PRC2}
	is 
	inspired by the conflict-driven clause learning (CDCL) 
	proposed in \cite{BS97,MS96,MS99} (see also \cite[Chapter 5]{BHMW09}).
CDCL
	enable us to find (or learn)
	a new pivot combination from the conflict evidence found by UP.

The performance of UP tends to increase 
	when 	
	combinations have good structures for finding conflict combinations under a good variable order:
	i.e., 
	not too few number of combinations,
	high frequency of small size combinations, 
	bias of occurrences of variables, 
	and so on.	
We have seen in \refTab{4_Tbl_Main_cal_bspZ}
	that 
	the runtimes of $\letMJPO$ and EDS 
	are 
	much better than those of $\letKNT$,
	which 
	seems to be due to the difference in numbers of relations (or combinations).
	
\end{remarkA}



\end{document}

%% file: letters.tex
	\newcommand{\letB}{\mathfrak{b}}
	\newcommand{\letD}{d}
	\newcommand{\letMJPO}{\mathrm{MJPO}}
	\newcommand{\letKNT}{\mathrm{KNT}}
	\newcommand{\dirM}{\mathtt{Main\_make}}	
	\newcommand{\dirC}{\mathtt{Main\_cal}}
	\newcommand{\dirD}{\mathtt{Data}}
	\newcommand{\dirW}{\mathtt{Work}}
	\newcommand{\fleM}[1][?]{\ifx0#1\mathtt{0\_preparation.py}\else
							\ifx1#1\mathtt{1\_product.py}\else
							\ifx2#1\mathtt{2\_regularized\_product.py}\else
							\ifx3#1\mathtt{3\_extended\_relation.py}\else
							\ifx4#1\mathtt{4\_binary\_system.py}\else !ERROR!
						\fi\fi\fi\fi\fi}
	\newcommand{\fleC}[1][d]{\ifx d#1\mathtt{dimensions.py}\else !ERROR! \fi}
	\newcommand{\fleConfig}[1][d]{\ifx d#1\mathtt{config.txt}\else !ERROR! \fi}